\documentclass[11pt,leqno]{amsart}
\usepackage{amssymb}

\newcommand{\C}{{\mathbb{C}}}

\newcommand{\N}{{\mathbb{N}}}
\newcommand{\Q}{{\mathbb{Q}}}

\newcommand{\T}{{\mathbb{T}}}
\newcommand{\Z}{{\mathbb{Z}}}

\newcommand{\cH}{{\mathcal{H}}}
\newcommand{\cL}{{\mathcal{L}}}

\newcommand{\Sym}{{\mathfrak{S}}}
\newcommand{\fC}{{\mathfrak{C}}}
\newcommand{\Irr}{{\operatorname{Irr}}}
\renewcommand{\leq}{\leqslant}
\renewcommand{\geq}{\geqslant}
\renewcommand{\atop}[2]{\genfrac{}{}{0pt}{}{#1}{#2}}

\def\infspe{\hspace{0.1em}\mathop{\preccurlyeq}\nolimits\hspace{0.1em}}

\mathchardef\inferieur="321E
\mathchardef\superieur="321F

\def\Ind{{\mathrm{Ind}}}

\def\lacril{~\longleftrightarrow_L~}
\def\lacrils#1{\stackrel{#1}{\longrightarrow}}
\def\Lacrils#1{\stackrel{#1}{\longleftrightarrow}}

\def\cedricl{~\leftrightsquigarrow_L~}

\def\finl{~$\SS \square$}
\def\matrice#1{\left(\begin{array}{ccccccccccccccccccc}#1\end{array}\right)}

\def\NM{{\mathbb{N}}}

\def\QM{{\mathbb{Q}}}

\def\ZM{{\mathbb{Z}}}

\def\LG{{\mathfrak L}}
\def\MG{{\mathfrak M}}

\def\SG{{\mathfrak S}}

\def\aG{{\mathfrak a}}

\def\a{\alpha}
\def\b{\beta}

\def\G{\Gamma}

\def\e{\varepsilon}

\def\ch{\chi}
\def\l{\lambda}

\def\m{\mu}

\def\s{\sigma}
\def\Sig{\Sigma}
\def\th{\theta}

\def\t{\tau}

\def\HC{{\mathcal{H}}}
\def\IC{{\mathcal{I}}}

\def\LC{{\mathcal{L}}}

\def\RC{{\mathcal{R}}}

\def\VC{{\mathcal{V}}}

\def\Mti{{\tilde{M}}}

\def\aha{{\hat{a}}}

\def\vha{{\hat{v}}}

\def\Vha{{\hat{V}}}

\def\thet{{\tilde{\theta}}}

\def\SS{\scriptstyle}

\def\Irr{\mathop{\mathrm{Irr}}\nolimits}
\def\Res{\mathop{\mathrm{Res}}\nolimits}

\makeatletter
\def\qedsymbol{\RIfM@\bgroup\else$\bgroup\aftergroup$\fi
  \vcenter{\hrule\hbox{\vrule\@height.6em\kern.6em\vrule}\hrule}\egroup}
\def\qed{\RIfM@\else\unskip\nobreak\fi\quad\qedsymbol}
\def\pf{%
\par\topsep6\p@\@plus6\p@ \trivlist
\item[\hskip\labelsep\it Proof.]\ignorespaces}

\makeatother

\newtheorem{thm}{Theorem}[section]
\newtheorem{lem}[thm]{Lemma}

\newtheorem{cor}[thm]{Corollary}
\newtheorem{prop}[thm]{Proposition}

\theoremstyle{definition}
\newtheorem{exmp}[thm]{Example}
\newtheorem{defn}[thm]{Definition}
\newtheorem{what}[thm]{}

\theoremstyle{remark}
\newtheorem{rem}[thm]{Remark}

\begin{document}

\title{Left cells in type $B_n$ with unequal parameters}

\author{C\'edric Bonnaf\'e and Lacrimioara Iancu}
\address{Laboratoire de Math\'ematiques de Besan\c{c}on,
Universit\'e de Franche-Comt\'e, 16 Route de Gray, 25030 Besan\c{c}on
Cedex, France} 

\makeatletter
\email{bonnafe@math.univ-fcomte.fr}
\address{Laboratoire de Math\'ematiques de Besan\c{c}on,
Universit\'e de Franche-Comt\'e, 16 Route de Gray, 25030 Besan\c{c}on
Cedex, France} 
\address{Facultatea de Stiinte, Universitatea de Nord Baia Mare, Victoriei 76,
RO-4800 Baia Mare, Romania} 

\makeatletter
\email{iancu@math.univ-fcomte.fr}
\makeatother


\date{}

\begin{abstract} Kazhdan and Lusztig have shown that the partition of
the symmetric group $\Sym_n$ into left cells is given by the
Robinson--Schensted correspondence. The aim of this paper is 
to provide a similar description of the left cells in type 
$B_n$ for a special class of choices of unequal parameters.
This is based on a generalization of the Robinson--Schensted 
correspondence in type $B_n$. We 
provide an explicit description of the left cell representations 
and show that they are irreducible and constructible.
\end{abstract}

\maketitle

\pagestyle{myheadings}

\markboth{Bonnaf\'e and Iancu}{Left cells in type $B_n$}

\section{Introduction} \label{sec1}

The Robinson--Schensted correspondence is a classical combinatorial
instrument which gives rise to a partition of the symmetric group
$\Sym_n$ into pieces which are indexed by the various standard tableaux 
of size $n$ (with a filling by the numbers $1,\ldots,n$). Kazhdan and Lusztig
have given a completely different description of that partition by using 
the construction of a new basis (the ``Kazhdan--Lusztig basis'') of 
the Iwahori--Hecke algebra of $\Sym_n$. In this context, the pieces in 
the partition are called {\em left cells}. Now the definition of left 
cells makes sense for any (finite or infinite) Coxeter group, using the 
Kazhdan--Lusztig basis of the one-parameter or even multi-parameter 
Iwahori--Hecke algebra. One of the important aspects of this construction
is that each left cell gives rise to a representation of the 
Iwahori--Hecke algebra where the underlying vector space has a natural 
basis indexed by the elements in that left cell. 

Now, in the case of finite Coxeter groups and one-parameter 
Iwahori--Hecke algebras, the decomposition of the left cell 
representations into irreducible representations is completely known. 
For the symmetric group $\Sym_n$, Kazhdan and Lusztig \cite{KaLu} 
showed that each left cell representation actually is irreducible. In 
the remaining types, the left cell representations are no longer
irreducible and Lusztig \cite{Lusztig86} showed how they decompose
into irreducibles. 

This paper is concerned with the multi-parameter case. Note that, as far 
as finite Coxeter groups are concerned, we only have to deal with 
the dihedral groups and Coxeter groups of type $F_4$ and $B_n$.
For type $B_n$ and a special choice of the parameters (which allows
a geometric interpretation), Lusztig \cite{Lusztig83} showed that, 
again, all left cell representations are irreducible. For the dihedral
groups and type $F_4$, results in the multi-parameter case have been 
obtained by Geck and Pfeiffer \cite{ourbuch}, Geck \cite{Ge1} and Lusztig
\cite{Lusztig02}.

In this paper, we consider Coxeter groups of type $B_n$ with diagram
and parameters of the corresponding Hecke algebra $\HC_n$ given as follows.

\begin{center}
\begin{picture}(300,30)
\put(  0, 10){$B_n$}
\put( 40, 10){\circle{10}}
\put( 44,  7){\line(1,0){33}}
\put( 44, 13){\line(1,0){33}}
\put( 81, 10){\circle{10}}
\put( 86, 10){\line(1,0){29}}
\put(120, 10){\circle{10}}
\put(125, 10){\line(1,0){20}}
\put(155,  7){$\cdot$}
\put(165,  7){$\cdot$}
\put(175,  7){$\cdot$}
\put(185, 10){\line(1,0){20}}
\put(210, 10){\circle{10}}
\put( 38, 22){$q^c$}
\put( 78, 22){$q$}
\put(118, 22){$q$}
\put(208, 22){$q$}
\put(250, 08){where $c\geq 1$}
\end{picture}
\end{center}

It is easy to see that the left cells are independent of the particular 
value of $c$, as long as $c$ is sufficiently large. We are precisely 
interested in this ``asymptotic'' case, where $c$ is sufficiently large. 
We shall actually prove a result which implies that all 
the left cell representations in that case are irreducible. This gives 
a new construction of ``integral'' forms of the irreducible representations; 
the first such construction is due to Dipper, James and Murphy 
\cite{DiJa95}. Our result involves a generalization of the classical
Robinson--Schensted correspondence.  This is the subject of 
Section~\ref{sec2}; the generalization does not only work for the 
Coxeter groups of type $B_n$ but for all complex reflection groups of 
type $G(e,1,n)$ (see \cite{iancu}).

Analogously to the case of $\Sym_n$, the generalized Robinson--Schensted 
correspondence gives rise to a partition of the Coxeter group $W_n$ of type $B_n$ 
into pieces which are indexed by pairs of standard bitableaux 
of total size $n$ (with a filling by the numbers $1,\ldots,n$).

Our aim will be to show that the left cells in type $B_n$ with the above
choice of the parameters are given by the generalized Robinson--Schensted 
correspondence. As a consequence, we obtain that the left cell representations 
of $\HC_n$ are irreducible and we retrieve the classical parametrisation 
of irreducible $\HC_n$-modules. 

The paper is organized as follows. After introducing the general set-up 
in \S\ref{setup}, we define in \S\ref{sec2} the generalized Robinson-Schensted 
correpondence and give its first properties (Knuth correspondence, 
compatibility with parabolic subgroups...). In \S \ref{sec3'}, we define 
a decomposition of elements of $W_n$ which sounds like {\it Clifford theory 
for elements}. This decomposition gives a new description of Robinson-Schensted 
cells in terms of cells for symmetric groups (Proposition \ref{equivalent}). 
One of the main tools developed in this section is Proposition \ref{bijections xnl}. 

In \S\ref{sec4}, we recall the basic notions and 
results concerning the construction of the Kazhdan--Lusztig basis and left
cells in the multi-parameter case for general Coxeter groups. For this purpose, 
it is convenient
to work in the general setting described by Lusztig 
(see \cite{Lusztig83} and \cite{Lusztig02}).
In \S\ref{sec5}, we come back to $W_n$: 
we will replace the parameter $q^c$ by a new variable $Q$ and
work with the Iwahori--Hecke algebra $\HC_n$ of type $B_n$ with two independent
parameters $q$ and $Q$. The main results of this section are the following.  
First, Kazhdan-Lusztig polynomials are polynomials only in $q$ (Theorem 
\ref{mythm} (a)). Secondly, we obtain a kind of grading 
for left cells (Theorems \ref{mythm} (b) and \ref{theo mt}). 
The last section is devoted to the proof of the main results 
of this paper, namely the explicit description of Kazhdan-Lusztig 
left cells (Theorem \ref{main theo}) and the fact that left cells 
representations are irreducible and explicitly determined 
(Proposition \ref{irr cell}).


\section{The set-up}\label{setup}

We introduce in this section all the notation we will need 
concerning the Coxeter group of type $B_n$. 
This group has a presentation with set of generators $S_n=\{ t,s_1,
\ldots, s_{n-1}\}$ and defining relations
\begin{equation*}
\left\{\begin{array}{r@{\hspace{1mm}}c@{\hspace{1mm}}l@{\hspace{1mm}}
c@{\hspace{1mm}}c@{\hspace{1mm}}l@{\hspace{1mm}}} t^2 & = & 1, & \quad 
s_i^2 & = & 1 \mbox{ for $i\geq 1$},\\t s_1t s_1 & = &s_1ts_1t,
& \quad  ts_i &=& s_it \mbox{ for $i \geq 2$},\\ s_is_{i+1}s_i & =& s_{i+1}
s_i s_{i+1} \mbox{ for $i \geq 1$}, &\quad s_is_j& =&s_js_i \mbox{ if $|i-j|
\geq 2$}.  \end{array} \right.
\end{equation*}
We visualize this presentation by the diagram
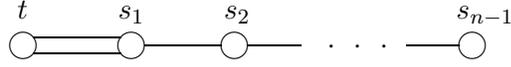
\begin{center}
\begin{picture}(220,30)
\put( 40, 10){\circle{10}}
\put( 44,  7){\line(1,0){33}}
\put( 44, 13){\line(1,0){33}}
\put( 81, 10){\circle{10}}
\put( 86, 10){\line(1,0){29}}
\put(120, 10){\circle{10}}
\put(125, 10){\line(1,0){20}}
\put(155,  7){$\cdot$}
\put(165,  7){$\cdot$}
\put(175,  7){$\cdot$}
\put(185, 10){\line(1,0){20}}
\put(210, 10){\circle{10}}
\put( 38, 20){$t$}
\put( 76, 20){$s_1$}
\put(116, 20){$s_2$}
\put(204, 20){$s_{n{-}1}$}
\end{picture}
\end{center}
A group with a presentation as above can be naturally realized as the
finite reflection group of type $G(2,1,n)$, that is, as the subgroup of 
$\text{GL}_n({\C})$ consisting of all matrices whose non-zero coefficients 
are $1$ or ${-}1$ and where there is precisely one non-zero 
coefficient in each row and each column. For our purposes, it will be 
more convenient to work with a different realisation, using permutations.

\begin{what} \label{abs22} {\bf $B_n$ as a permutation 
group.} Let $n\geq 1$ and consider the set 
$$I_n=I_n^+\cup I_n^-,$$
$$I_n^+=\{ 1,2,\ldots ,n\}\quad\text{and}\quad
I_n^-=-I_n^+=\{ {-}1,{-}2,\ldots ,{-}n\}.\leqno{\text{where}}$$
We denote by $\SG(I_n)$ the group of permutations of the set $I_n$ and 
we set 
$$W_n:=\{\pi \in \SG(I_n)~|~\forall i \in I_n,~\pi(-i)=-\pi(i)\}.$$
In other words, if $w_n \in \SG(I_n)$ is defined by $I_n \to I_n$, $i \mapsto -i$, 
then $W_n$ is the centralizer of $w_n$ in $\SG(I_n)$. 
We define the following transpositions in $\SG(I_n)$:
\begin{equation*}
t_i:=(i,\, {-}i) \quad \mbox{for $1\leq 
i\leq n$}.
\end{equation*}
Then $w_n=t_1t_2\cdots t_n$. 
The order formulas for centralizers in symmetric groups 
(see \cite[Chap.~I]{Mac}) show that $|W_n|=2^n\cdot n!$. It is easily
checked that $t_1,t_2,\ldots,t_n$ generate a subgroup $N_n\subseteq 
W_n$ which is isomorphic to $(\ZM/2\ZM)^n$. Furthermore, the elements
\begin{equation*}
\left\{\begin{array}{rcl}
s_1 & := &(1,2) \cdot ({-}1, {-}2),\\
s_2 & := &(2,3) \cdot ({-}2, {-}3),\\
& \vdots &\\
s_{n-1} & := & ((n-1),n) \cdot ({-}(n-1), {-}n),
\end{array}\right.
\end{equation*}
generate a subgroup $\SG_n \subseteq W_n$ which is isomorphic to 
the symmetric group of degree $n$. We set $\Sig_n=\{s_1,s_2,\cdots,s_{n-1}\}$. 
Since $N_n\cap \SG_n=\{1\}$, we have $|\SG_n N_n|=|\SG_n|\cdot |N_n|=2^n\cdot n!$ 
and so $W_n=N_n\SG_n$. We now set
$$S_n=\{t,s_1,s_2,\cdots,s_{n{-}1}\}\quad\text{where}~t=t_1.$$
The previous discussion shows that 
\[W_n=\langle S_n \rangle.\]
In terms of these generators, the above transpositions $t_i$ are given 
recursively by $t_{i+1}=s_it_is_i$ for $1\leq i \leq n-1$.
Finally, since the generators $t,s_1,\ldots,s_{n-1}$ satisfy the relations 
specified by the above diagram and since $|W_n|=2^n\cdot n!$, we 
conclude that these relations form a set of defining relations for $W_n$.
\end{what} 

\begin{rem} \label{nicearray} The fact that each element $\pi \in W_n$ 
commutes with $w_n=t_1t_2\cdots t_n$ implies that $\pi$ is uniquely 
determined by the images of $1,\ldots,n$. Indeed, if we know $\pi(i)$, 
then we also know $\pi(-i)=-\pi(i)$.
Thus, if we set
\[ \pi(i)=\varepsilon_i \cdot p_i \qquad \mbox{where $1\leq i \leq n$ and
$\varepsilon_i \in\{ {+}1, {-}1\}$},\]
then the sequence of numbers $p_1,\ldots,p_n$ is a permutation of $1,\ldots,
n$. Thus, we shall represent $\pi$ by the array  
\[\pi =\begin{pmatrix} 1 & 2 & \ldots & n\\ \varepsilon_1 \cdot p_1& 
\varepsilon_2 \cdot p_2  & \ldots & \varepsilon_n \cdot p_n \end{pmatrix} 
\in W_n.~{\scriptstyle{\square}}\] 
\end{rem}

\begin{what}\label{subsection length}{\bf Length function, Bruhat ordering.} 
Using this choice of generators, we can define the length function 
$\ell : W_n \to \NM$. It is easily checked that $w_n$ is the longest 
element of $W_n$~: we have $\ell(w_n)=n^2$. If $w \in W_n$, we denote 
by $\ell_t(w)$ the number of occurences of $t$ in a reduced decomposition 
of $w$. This does not depend on the choice of the reduced decomposition. 
We set $\ell_s(w)=\ell(w)-\ell_t(w)$. It is easily checked that, if $w$ and $w'$ 
are elements of $W_n$ such that $\ell(ww')=\ell(w)+\ell(w')$, then 
$\ell_t(ww')=\ell_t(w)+\ell_t(w')$ and $\ell_s(ww')=\ell_s(w)+\ell_s(w')$. 
For instance, we have, for every $i \in I_n$, 
$$\ell_t(t_i)=1, \quad \ell_s(t_i)=2(i-1) \quad \text{and}\quad \ell(t_i)=2i-1.$$
So $\ell_t(w_n)=n$ and $\ell_s(w_n)=n^2-n$. This implies that, for every $w \in W_n$, 
we have $\ell_t(w) \leq n$. 

We denote by $\leq$ the Bruhat order on $W_n$ defined by the set of generators 
$S_n$. We write $x < y$ to say that $x \leq y$ and $x \not= y$. 
If $w \in W_n$, we define its left descent set $\LC(w)$ and its right descent 
set $\RC(w)$ as follows~:
$$\LC(w)=\{s \in S_n~|~sw < w\}$$
$$\RC(w)=\{s \in S_n~|~ws < w\}.\leqno{\mathrm{and}}$$
\end{what}

\section{On the Robinson--Schensted correspondence} \label{sec2}

In this section we describe a generalization of the classical 
Robinson--Schensted correspondence (which is concerned with
the symmetric group $\SG_n$) to the Coxeter group of type $B_n$.
For more details on the classical correspondence see Knuth 
\cite[5.1.4]{Knuth3} or Fulton \cite[Part~I]{Fulton}. 
\begin{what} {\bf A generalized Robinson-Schensted correspondence.}
Let us first introduce some more notation. If 
$\l$ is a partition, and if $T$ is a standard tableau of 
shape $\l$, we set $|T|=|\l|$ (the number $|\l|$ is called the {\it size} 
of $T$). A {\it bipartition $($of $n$)} is a pair 
$(\l,\m)$ of partitions (such that $|\l|+|\m|=n$). 
A {\it bitableau} is a pair of tableaux. If $(T_1,T_2)$ is a bitableau 
such that $T_1$ is of shape $\l$ and $T_2$ is of shape $\m$, we say 
that $(\l,\m)$ is the {\it shape} of $(T_1,T_2)$ and that 
$|\l|+|\m|$ is the size of $(T_1,T_2)$. The bitableau is said to be 
{\it $n$-standard} if $T_1$ and $T_2$ are standard tableaux, 
if $|T_1|+|T_2|=n$ and if the filling of $T_1$ and $T_2$ is 
the set $I_n^+$.

In order to generalize the Robinson-Schensted correspondence to 
$W_n$, we work with the realisation of $W_n$ as a subgroup 
of $\SG(I_n)$ and use Remark~\ref{nicearray} to represent the elements 
of $W_n$. Thus, let $\pi\in W_n$. 
Then we define a pair of $n$-standard bitableaux:
\[ \renewcommand{\arraystretch}{1.5}
(A_n(\pi),B_n(\pi)) \quad \mbox{where}\quad
\left\{\begin{array}{ll} 
A_n(\pi)=(A_n^{+}(\pi ),\,A_n^{-}(\pi ))\\
B_n(\pi)=(B_n^{+}(\pi ),\,B_n^{-}(\pi ))
\end{array}\right.\]
and $A_n(\pi)$, $B_n(\pi)$ have the same shape. This is done as follows. 

Apply the Knuth insertion procedure as follows: insert succesively
the numbers $p_i$ into two initially empty tableaux $A_n^+ (\pi ),
\, A_n^- (\pi )$, more precisely insert $p_i$ into $A_n^{\varepsilon _i}
(\pi )$.
Note that at each step this yields a new box, located on the $a_i$th
row and $b_i$th column, say, of $A_n^{\varepsilon _i}(\pi )$. Now add a 
box containing $i$ to $B_n^{\varepsilon _i}(\pi )$ on its $(a_i,b_i)$ 
position ("keep the record").
 

\begin{exmp}\label{ilx en faut bien un}
Let us consider an element $\pi\in W_7$ represented, as in 
Remark \ref{nicearray}, by the array
$$\pi=\matrice{1 & 2 & 3 & 4 & 5 & 6 & 7 \\ -4 & 3 & 6 & -1 & 
7 & -2 & 5}.$$ Then $(A(\pi),B(\pi))$ is equal to 
\[
\begin{picture}(330,40)
\put( 0, 15){$A(\pi ) =$}
\put( 40,36){\line(1,0){54}}
\put( 40, 18){\line(1,0){54}} 
\put( 40, 36){\line(0,-1){36}}
\put(58, 36){\line(0,-1){36}}
\put(76, 18){\line(0,1){18}}
\put(94, 18){\line(0,1){18}}
\put(40, 0){\line(1,0){18}}
\put(46, 23){3}
\put(64, 23){5} 
\put(82, 23){7}
\put(46, 5){6}

\put( 110,36){\line(1,0){36}}
\put( 110, 18){\line(1,0){36}} 
\put( 110, 36){\line(0,-1){36}}
\put(128, 36){\line(0,-1){36}}
\put(146, 18){\line(0,1){18}}
\put(110, 0){\line(1,0){18}}
\put(116, 23){1}
\put(134, 23){2}
\put(116, 5){4}

\put(180, 15){$B(\pi ) =$}
\put( 220,36){\line(1,0){54}}
\put( 220, 18){\line(1,0){54}} 
\put( 220, 36){\line(0,-1){36}}
\put(238, 36){\line(0,-1){36}}
\put(256, 18){\line(0,1){18}}
\put(274, 18){\line(0,1){18}}
\put(220, 0){\line(1,0){18}}
\put( 226, 23){2}
\put(244, 23){3}
\put(262, 23){5}
\put(226, 5){7}

\put( 290,36){\line(1,0){36}}
\put( 290, 18){\line(1,0){36}} 
\put( 290, 36){\line(0,-1){36}}
\put(308, 36){\line(0,-1){36}}
\put(326, 18){\line(0,1){18}}
\put(290, 0){\line(1,0){18}}
\put(296, 23){1}
\put(314, 23){6}
\put(296, 5){4}
\end{picture}
\]
\finl
 \end{exmp}

\begin{thm}\label{rs1} With the above notation, the following hold.
\begin{itemize}
\item[(a)] The map $\pi \mapsto (A_n(\pi),B_n(\pi))$ is 
a bijection from $W_n$ onto the set of all pairs of
$n$-standard bitableaux of the same shape. 
\item[(b)] For any $\pi \in W_n$, we have $A_n(\pi^{-1})=
B_n(\pi)$ and $B_n(\pi^{-1})=A_n(\pi)$.
\item[(c)] The number of generalized Robinson--Schensted cells equals the 
number of involutions in $W_n$.
\end{itemize}
\end{thm}

\begin{proof} The proof uses essentially the same argument as for 
the symmetric group \cite[5.1.4, Theorems A and B]{Knuth3}. 
(For a description of other possible generalizations see 
\cite[4.2.3]{vLee}.)
(a) It is clear from the construction that, for any $\pi \in W_n$, $A_n(\pi)$ and 
$B_n(\pi)$ always have the same shape; furthermore: (i) 
since they are obtained through the bumping procedure, $A_n^{+}(\pi)$ 
and $A_n^{-}(\pi)$ are standard tableaux, (ii) $B_n^{+}(\pi)$ and $B_n^{-}(\pi)$ are 
standard tableaux too, since we always add elements on their periphery
in increasing order. 

Conversely, given a pair $(A,B)$ of $n$-standard bitableaux of the same shape 
(with $A=(A^+,A^-)$ and $B=(B^+,B^-)$), we can find the corresponding
array (and so the element $\pi\in W_{n}$) as follows: for 
$i=n,\ldots ,2,1$, let $(\nu ,a,b)$ be defined by the fact that the number
$i$ appears in the $B^\nu$ tableau, on its $a$th row and $b$th column. 
We set then $\varepsilon_{i}=\nu$. Now let $p_{i}$ be the element $x$ 
that is removed when applying the deleting algorithm (inverse of the 
insertion algorithm) to the $(\nu ,a,b)$ box of $A$. 
The two constructions we have described are inverses of each other.

(b) The argument is identical to the one in \cite[5.1.4, Theorem B]{Knuth3}.

(c) is an easy consequence of (a) and (b).
\end{proof}

\begin{defn} \label{genRS} Let $\T$ be an $n$-standard bitableau.
Then the set
\[ \T (W_n):=\{w\in W_n\mid B_n(w)=\T\}\]
is non-empty and will be called a {\em generalized Robinson--Schensted 
cell} (or generalized RS-cell for short) of $W_n$. We have a partition
\[W_n=\coprod_{\T} \T(W_n),\]  where $\T$ runs over the set of 
all $n$-standard bitableaux. An example is given in Table~\ref{tab1}.
\end{defn}
\end{what}

\begin{rem}\label{t longueur}
If $\pi \in W_n$ then, by construction, $\ell_t(\pi)=|A_n^-(\pi)|=
|B_n^-(\pi)|$. Therefore, the function $\ell_t$ is constant on 
generalized Robinson-Schensted cells. Moreover, if $\pi \in \SG_n$, 
then $A_n^-(\pi)=B_n^-(\pi)=\varnothing$ and $(A_n^+(\pi),
B_n^+(\pi))$ is the pair of standard tableaux associated 
to $\pi$ via the usual Robinson-Schensted correspondence.\finl
\end{rem}

\begin{exmp}\label{tableau w0}
We have 
\[ 
\begin{picture}(330,25)
\put(0, 5){$A_n (1) =$}
\put(45, 18){\line(1,0){80}}
\put(45, 0){\line(1,0){80}} 
\put(45, 0){\line(0,1){18}}
\put(63, 0){\line(0,1){18}}
\put(81, 0){\line(0,1){18}}
\put(107, 0){\line(0,1){18}}
\put(125, 0){\line(0,1){18}}
\put(51, 5){1}
\put(69, 5){2}
\put(87, 5){...}
\put(113, 5){n}

\put( 140,7){$\varnothing $}

\put(180, 5){$B_n(1) =$}
\put(225,18){\line(1,0){80}}
\put(225, 0){\line(1,0){80}}
\put(225, 0){\line(0,1){18}}
\put(243, 0){\line(0,1){18}}
\put(261, 0){\line(0,1){18}}
\put(287, 0){\line(0,1){18}}
\put(305, 0){\line(0,1){18}}
\put(231, 5){1}
\put(249, 5){2}
\put(267, 5){...}
\put(293, 5){n}

\put(320,7){$\varnothing $}
\end{picture}
\]
On the other hand, 
\[ 
\begin{picture}(350,25) 
\put(0, 5){$A_n (w_n)= $}
\put(55, 5){$\varnothing $}

\put(80, 18){\line(1,0){80}}
\put(80, 0){\line(1,0){80}}
\put(80, 0){\line(0,1){18}}
\put(98, 0){\line(0,1){18}}
\put(116, 0){\line(0,1){18}}
\put(142, 0){\line(0,1){18}}
\put(160, 0){\line(0,1){18}}
\put(86, 5){1}
\put(104, 5){2}
\put(122, 5){...}
\put(148, 5){n}

\put(190, 5){$B_n (w_n)=$}
\put(245, 5){$\varnothing $}

\put(270,18){\line(1,0){80}}
\put(270, 0){\line(1,0){80}}
\put(270, 0){\line(0,1){18}}
\put(288, 0){\line(0,1){18}}
\put(306, 0){\line(0,1){18}}
\put(332, 0){\line(0,1){18}}
\put(350, 0){\line(0,1){18}}
\put(276, 5){1}
\put(294, 5){2}
\put(312, 5){...}
\put(338, 5){n}

\end{picture}
\]

Moreover, $\{1\}$ and $\{w_n\}$ are 
generalized Robinson-Schensted left cells.\finl
\end{exmp}

\begin{table}[htbp]  \caption{Generalized Robinson--Schensted cells in 
type $B_3$} \label{tab1}
\begin{center} \unitlength 0.55pt
$\begin{array}{l@{\hspace*{1cm}}rl} 
\hline\mbox{Generalized RS-cell} 
&\multicolumn{2}{l}{\mbox{$B$-bitableaux}}
\\\hline  \{1\}  &   
\begin{picture}(70,25)
\put( 0,20){\line(1,0){60}}
\put( 0, 0){\line(1,0){60}}
\put( 0, 0){\line(0,1){20}}
\put(20, 0){\line(0,1){20}}
\put(40, 0){\line(0,1){20}}
\put(60, 0){\line(0,1){20}}
\put( 5, 5){1}
\put(25, 5){2}
\put(45, 5){3}
\end{picture}  & 
\begin{picture}(20,25) \put(0,7){\line(1,0){20}} \end{picture}
\\ \hline 
 \{ s_{2},\quad s_{1}s_{2}\}  &  
\begin{picture}(50,50)
\put( 0,20){\line(1,0){40}}
\put( 0,40){\line(1,0){40}}
\put( 0, 0){\line(1,0){20}}
\put( 0, 0){\line(0,1){40}}
\put(20, 0){\line(0,1){40}}
\put(40,20){\line(0,1){20}}
\put( 5,25){1}
\put(25,25){2}
\put( 5, 5){3}
\end{picture}  & 
\begin{picture}(20,25) \put(0,17){\line(1,0){20}} \end{picture}
\\ \{s_{1},\quad s_{2}s_{1}\}  &  
\begin{picture}(50,45)
\put( 0,20){\line(1,0){40}}
\put( 0,40){\line(1,0){40}}
\put( 0, 0){\line(1,0){20}}
\put( 0, 0){\line(0,1){40}}
\put(20, 0){\line(0,1){40}}
\put(40,20){\line(0,1){20}}
\put( 5,25){1}
\put(25,25){3}
\put( 5, 5){2}
\end{picture}  & 
\begin{picture}(20,25) \put(0,17){\line(1,0){20}} \end{picture}
\\ \hline
  \{s_{1}s_{2}s_{1}\}  &  
\begin{picture}(30,65)
\put( 0,60){\line(1,0){20}}
\put( 0,40){\line(1,0){20}}
\put( 0,20){\line(1,0){20}}
\put( 0, 0){\line(1,0){20}}
\put( 0, 0){\line(0,1){60}}
\put(20, 0){\line(0,1){60}}
\put( 5, 5){3}
\put( 5,25){2}
\put( 5,45){1}
\end{picture}  & 
\begin{picture}(20,25) \put(0,27){\line(1,0){20}} \end{picture}
\\ \hline
 \{ t,\quad  s_{1}t,\quad  s_{2}s_{1}t\}  &  
\begin{picture}(50,25)
\put( 0,20){\line(1,0){40}}
\put( 0, 0){\line(1,0){40}}
\put( 0, 0){\line(0,1){20}}
\put(20, 0){\line(0,1){20}}
\put(40, 0){\line(0,1){20}}
\put( 5, 5){2}
\put(25, 5){3}
\end{picture}  & 
\begin{picture}(20,25) 
\put( 0,20){\line(1,0){20}}
\put( 0, 0){\line(1,0){20}}
\put( 0, 0){\line(0,1){20}}
\put(20, 0){\line(0,1){20}}
\put( 5, 5){1}
\end{picture}
\\ \{ts_{1},\quad s_{1}ts_{1},\quad s_{2}s_{1}ts_{1}\}  &  
\begin{picture}(50,25)
\put( 0,20){\line(1,0){40}}
\put( 0, 0){\line(1,0){40}}
\put( 0, 0){\line(0,1){20}}
\put(20, 0){\line(0,1){20}}
\put(40, 0){\line(0,1){20}}
\put( 5, 5){1}
\put(25, 5){3}
\end{picture}  & 
\begin{picture}(20,25) 
\put( 0,20){\line(1,0){20}}
\put( 0, 0){\line(1,0){20}}
\put( 0, 0){\line(0,1){20}}
\put(20, 0){\line(0,1){20}}
\put( 5, 5){2}
\end{picture}
\\ \{ts_{1}s_{2},\quad s_{1}ts_{1}s_{2},\quad s_{2}s_{1}ts_{1}s_{2}\}  &  
\begin{picture}(50,25)
\put( 0,20){\line(1,0){40}}
\put( 0, 0){\line(1,0){40}}
\put( 0, 0){\line(0,1){20}}
\put(20, 0){\line(0,1){20}}
\put(40, 0){\line(0,1){20}}
\put( 5, 5){1}
\put(25, 5){2}
\end{picture}  & 
\begin{picture}(20,25) 
\put( 0,20){\line(1,0){20}}
\put( 0, 0){\line(1,0){20}}
\put( 0, 0){\line(0,1){20}}
\put(20, 0){\line(0,1){20}}
\put( 5, 5){3}
\end{picture}
\\ \hline
  \{ts_{2},\quad s_{1}ts_{2},\quad s_{1}s_{2}s_{1}t\}  &  
\begin{picture}(30,45)
\put( 0,40){\line(1,0){20}}
\put( 0,20){\line(1,0){20}}
\put( 0, 0){\line(1,0){20}}
\put( 0, 0){\line(0,1){40}}
\put(20, 0){\line(0,1){40}}
\put( 5,25){2}
\put( 5, 5){3}
\end{picture}  & 
\begin{picture}(20,45) 
\put( 0,40){\line(1,0){20}}
\put( 0,20){\line(1,0){20}}
\put( 0,20){\line(0,1){20}}
\put(20,20){\line(0,1){20}}
\put( 5,25){1}
\end{picture}
\\ \{ts_{2}s_{1},\quad s_{1}ts_{2}s_{1},\quad s_{1}s_{2}s_{1}ts_{1}\}  &  
\begin{picture}(30,45)
\put( 0,40){\line(1,0){20}}
\put( 0,20){\line(1,0){20}}
\put( 0, 0){\line(1,0){20}}
\put( 0, 0){\line(0,1){40}}
\put(20, 0){\line(0,1){40}}
\put( 5,25){1}
\put( 5, 5){3}
\end{picture}  & 
\begin{picture}(20,45) 
\put( 0,40){\line(1,0){20}}
\put( 0,20){\line(1,0){20}}
\put( 0,20){\line(0,1){20}}
\put(20,20){\line(0,1){20}}
\put( 5,25){2}
\end{picture}
\\ \{ts_{1}s_{2}s_{1},\quad s_{1}ts_{1}s_{2}s_{1},
\quad s_{1}s_{2}s_{1}ts_{1}s_{2}\} &  
\begin{picture}(30,45)
\put( 0,40){\line(1,0){20}}
\put( 0,20){\line(1,0){20}}
\put( 0, 0){\line(1,0){20}}
\put( 0, 0){\line(0,1){40}}
\put(20, 0){\line(0,1){40}}
\put( 5,25){1}
\put( 5, 5){2}
\end{picture}  & 
\begin{picture}(20,45) 
\put( 0,40){\line(1,0){20}}
\put( 0,20){\line(1,0){20}}
\put( 0,20){\line(0,1){20}}
\put(20,20){\line(0,1){20}}
\put( 5,25){3}
\end{picture}
\\\hline \{ts_{1}ts_{1}s_{2}s_{1},\quad ts_{1}s_{2}s_{1}ts_{1}s_{2}, \quad 
s_{1}ts_{1}s_{2}s_{1}ts_{1}s_{2} \} &  
\begin{picture}(30,25) 
\put( 0,20){\line(1,0){20}}
\put( 0, 0){\line(1,0){20}}
\put( 0, 0){\line(0,1){20}}
\put(20, 0){\line(0,1){20}}
\put( 5, 5){1}
\end{picture}& 
\begin{picture}(40,25)
\put( 0,20){\line(1,0){40}}
\put( 0, 0){\line(1,0){40}}
\put( 0, 0){\line(0,1){20}}
\put(20, 0){\line(0,1){20}}
\put(40, 0){\line(0,1){20}}
\put( 5, 5){2}
\put(25, 5){3}
\end{picture}  
\\\{ts_{1}ts_{1}s_{2},\quad ts_{2}s_{1}ts_{1}s_{2},\quad 
s_{1}ts_{2}s_{1}ts_{1}s_{2}\}  &  
\begin{picture}(30,25) 
\put( 0,20){\line(1,0){20}}
\put( 0, 0){\line(1,0){20}}
\put( 0, 0){\line(0,1){20}}
\put(20, 0){\line(0,1){20}}
\put( 5, 5){2}
\end{picture}& 
\begin{picture}(40,25)
\put( 0,20){\line(1,0){40}}
\put( 0, 0){\line(1,0){40}}
\put( 0, 0){\line(0,1){20}}
\put(20, 0){\line(0,1){20}}
\put(40, 0){\line(0,1){20}}
\put( 5, 5){1}
\put(25, 5){3}
\end{picture}  
\\ \{ts_{1}ts_{1},\quad ts_{2}s_{1}ts_{1},\quad s_{1}ts_{2}s_{1}ts_{1}\}&  
\begin{picture}(30,25) 
\put( 0,20){\line(1,0){20}}
\put( 0, 0){\line(1,0){20}}
\put( 0, 0){\line(0,1){20}}
\put(20, 0){\line(0,1){20}}
\put( 5, 5){3}
\end{picture}& 
\begin{picture}(40,25)
\put( 0,20){\line(1,0){40}}
\put( 0, 0){\line(1,0){40}}
\put( 0, 0){\line(0,1){20}}
\put(20, 0){\line(0,1){20}}
\put(40, 0){\line(0,1){20}}
\put( 5, 5){1}
\put(25, 5){2}
\end{picture}  
\\ \hline
\{ts_{1}ts_{2}s_{1},\quad ts_{1}s_{2}s_{1}ts_{1},\quad
s_{1}ts_{1}s_{2}s_{1}ts_{1}\}&
\begin{picture}(30,45) 
\put( 0,40){\line(1,0){20}}
\put( 0,20){\line(1,0){20}}
\put( 0,20){\line(0,1){20}}
\put(20,20){\line(0,1){20}}
\put( 5,25){1}
\end{picture}  & 
\begin{picture}(30,45)
\put( 0,40){\line(1,0){20}}
\put( 0,20){\line(1,0){20}}
\put( 0, 0){\line(1,0){20}}
\put( 0, 0){\line(0,1){40}}
\put(20, 0){\line(0,1){40}}
\put( 5,25){2}
\put( 5, 5){3}
\end{picture}
\\ \{ts_{1}ts_{2},\quad ts_{1}s_{2}s_{1}t,\quad s_{1}ts_{1}s_{2}s_{1}t\}&  
\begin{picture}(30,45) 
\put( 0,40){\line(1,0){20}}
\put( 0,20){\line(1,0){20}}
\put( 0,20){\line(0,1){20}}
\put(20,20){\line(0,1){20}}
\put( 5,25){2}
\end{picture}  & 
\begin{picture}(30,45)
\put( 0,40){\line(1,0){20}}
\put( 0,20){\line(1,0){20}}
\put( 0, 0){\line(1,0){20}}
\put( 0, 0){\line(0,1){40}}
\put(20, 0){\line(0,1){40}}
\put( 5,25){1}
\put( 5, 5){3}
\end{picture}
\\\{ts_{1}t,\quad  ts_{2}s_{1}t,\quad  s_{1}ts_{2}s_{1}t\}  &  
\begin{picture}(30,45) 
\put( 0,40){\line(1,0){20}}
\put( 0,20){\line(1,0){20}}
\put( 0,20){\line(0,1){20}}
\put(20,20){\line(0,1){20}}
\put( 5,25){3}
\end{picture}  & 
\begin{picture}(30,45)
\put( 0,40){\line(1,0){20}}
\put( 0,20){\line(1,0){20}}
\put( 0, 0){\line(1,0){20}}
\put( 0, 0){\line(0,1){40}}
\put(20, 0){\line(0,1){40}}
\put( 5,25){1}
\put( 5, 5){2}
\end{picture}
\\\hline \{ts_{1}ts_{2}s_{1}t\}  &  
\begin{picture}(30,65) \put(0,27){\line(1,0){20}} \end{picture} & 
\begin{picture}(30,65)
\put( 0,60){\line(1,0){20}}
\put( 0,40){\line(1,0){20}}
\put( 0,20){\line(1,0){20}}
\put( 0, 0){\line(1,0){20}}
\put( 0, 0){\line(0,1){60}}
\put(20, 0){\line(0,1){60}}
\put( 5, 5){3}
\put( 5,25){2}
\put( 5,45){1}
\end{picture}  
\\\hline \{ts_{1}ts_{1}s_{2}s_{1}t,\quad ts_{1}ts_{2}s_{1}ts_{1}s_{2}\} &  
\begin{picture}(30,45) \put(0,17){\line(1,0){20}} \end{picture} & 
\begin{picture}(40,45)
\put( 0,20){\line(1,0){40}}
\put( 0,40){\line(1,0){40}}
\put( 0, 0){\line(1,0){20}}
\put( 0, 0){\line(0,1){40}}
\put(20, 0){\line(0,1){40}}
\put(40,20){\line(0,1){20}}
\put( 5,25){1}
\put(25,25){3}
\put( 5, 5){2}
\end{picture}  
\\ \{ts_{1}ts_{2}s_{1}ts_{1},\quad ts_{1}ts_{1}s_{2}s_{1}ts_{1}\}  &  
\begin{picture}(30,45) \put(0,17){\line(1,0){20}} \end{picture} & 
\begin{picture}(40,45)
\put( 0,20){\line(1,0){40}}
\put( 0,40){\line(1,0){40}}
\put( 0, 0){\line(1,0){20}}
\put( 0, 0){\line(0,1){40}}
\put(20, 0){\line(0,1){40}}
\put(40,20){\line(0,1){20}}
\put( 5,25){1}
\put(25,25){2}
\put( 5, 5){3}
\end{picture}  
\\ \hline \{ts_{1}ts_{1}s_{2}s_{1}ts_{1}s_{2}\}  &  
\begin{picture}(30,25) \put(0,7){\line(1,0){20}} \end{picture} & 
\begin{picture}(60,25)
\put( 0,20){\line(1,0){60}}
\put( 0, 0){\line(1,0){60}}
\put( 0, 0){\line(0,1){20}}
\put(20, 0){\line(0,1){20}}
\put(40, 0){\line(0,1){20}}
\put(60, 0){\line(0,1){20}}
\put( 5, 5){1}
\put(25, 5){2}
\put(45, 5){3}
\end{picture} \\ \hline \end{array}$
\end{center}
\end{table} 

\begin{what} {\bf A generalization of a theorem of Knuth.} \label{sec3}
Knuth has given a purely group theoretical description of the partition
of $\Sym_n$ into Robinson--Schensted cells.  We wish to 
generalize that statement to the Coxeter groups of type $B_n$. 
We begin with some general definitions. First, we set 
$S_n^\prime=S_n \cup \{t_1,t_2,\cdots,t_n\}$. Then the extended 
left descent set of $w \in W_n$ is defined by 
$$\cL'(w):=\{ u\in S_n^\prime~|~ \ell(uw)<\ell(w)\}=
\{ u\in S_n^\prime\mid uw<w\}.$$ 
Let $x,y\in W_n$ and $s\in \Sig_n=S_n-\{t\}$; 
then we define $$x\lacrils{s}_L
y\qquad\stackrel{\text{def}}{\Longleftrightarrow} \qquad y=sx,\:  
\ell (y)>\ell (x) \mbox{ and } \cL'(x)\nsubseteq\cL'(y)\, ,$$ 
and we write $x\Lacrils{s}_L y$ if 
$x\lacrils{s}_L y$ or  $y\lacrils{s}_L x$.
Finally, we write $x\lacril y$ if there exists a sequence
$x=x_1,x_2,\ldots ,x_k =y$ and $s_i\in \Sig_n$ such that
$x_i\Lacrils{s_i}_L x_{i{+}1}$ for all $i$.

\begin{prop} \label{caracRS}
Let $x,y\in W_n$. Then \[ B_n(x)= B_n(y) \quad 
\mbox{ if and only if } \quad x\lacril y\, .\]
\end{prop}

\begin{proof} We first define "admissible transformations" in $W_n$.
Let $x\in W_n$ be represented, as in \ref{nicearray}, by the array
\begin{equation*}
x =\begin{pmatrix} 1 & 2 & \ldots & i & i{+}1 & \ldots & n\\
\varepsilon_1\cdot p_{1} & \varepsilon_2\cdot p_{2} & \ldots 
& \varepsilon_i\cdot p_i & \varepsilon_{i{+}1}\cdot p_{i{+}1} 
& \ldots & \varepsilon_n\cdot p_{n}
\end{pmatrix} 
\end{equation*}
with $p_j\in I_n^+$ and $\varepsilon_j\in\{ +1,-1\}$; interchanging 
$\varepsilon_i p_i$ and $\varepsilon_{i{+}1} p_{i{+}1}$ is an 
admissible transformation if we are in one of the following 
situations:

(a) $2\leq i\leq n{-}1,\; \varepsilon_{i{-}1}=\varepsilon_i =
\varepsilon_{i{+}1}$ and $p_{i{-}1}$ lies between $p_i$ and 
$p_{i{+}1}$ 

(b) $1\leq i\leq n{-}2,\; \varepsilon_i =\varepsilon_{i{+}1}= 
\varepsilon_{i{+}2} $ and $p_{i{+}2}$ lies between $p_i$ and    
$p_{i{+}1}$ 

(c) $\varepsilon_i=-\varepsilon_{i{+}1}$.

Now note that (c)-type transformations do not change the relative ordering 
of the numbers $x(i)$ belonging to $I_n^{+}$ or to $I_n^{-}$ respectively. 
Then, by applying \cite[5.1.4, Ex.~4]{Knuth3} we get
that two elements $x,y\in W_n$ have the same
$A$-bitableaux if and only if each of them can be obtained from the 
other through a finite number of admissible transformations. 

As for the group theoretical description of these admissible 
transformations, we obtain the following as an easy consequence of 
\cite[Ex. 9.10]{Taylor}:
$$\begin{array}{lcll}
\bullet\quad xs_{i} < x &\Longleftrightarrow &&  x(i),\,
x(i{+}1)\in I_n^{+} \mbox{ and } x(i{+}1)< x(i)\\
&&\mbox{ or } & x(i),\, x(i{+}1)\in I_n^{-} \mbox{ and } x(i{+}1)<x(i)\\
&&\mbox{ or } & x(i)\in I_n^{+},\; x(i{+}1)\in I_n^{-},\\
\bullet\quad \ell (xt_{i}) < \ell (x) &\Longleftrightarrow && x(i)\in I_n^{-}.
\end{array}$$
The proof of the proposition is now complete.
\end{proof}
\end{what}

\section{A third construction of generalized Robinson-Schensted cells}\label{sec3'}

\bigskip

\begin{what}\label{preliminaires abc}{\bf Preliminaries.} 
The parabolic subgroup $\Sym_n$ of $W_n$ generated by $\Sigma_n=
\{ s_1,\ldots ,s_{n{-}1}\}$ is isomorphic to the symmetric group
of degree $n$. We denote by $X_n$ the set of elements $w \in W_n$ 
which are of minimal length in $w\SG_n$ (they are usually
called {\em distinguished left coset representatives} of $\Sym_n$ in
$W_n$). If $0 \leq l \leq n$, we set
$$X_n^{(l)}=\{w \in X_n~|~\ell_t(w)=l\}.$$
Let us give a description of $X_n^{(l)}$. We define $r_1=t$ and, for 
$1 \leq i \leq n-1$, we set $r_{i+1} = s_i r_i=s_i \dots s_2 s_1 t$. 
Then
$$X_n^{(l)}=\{r_{i_1} r_{i_2} \dots r_{i_l}~|~1 \leq i_1 < i_2 < \dots 
< i_l \leq n\}.$$
We have for instance $X_n^{(0)}=\{1\}$. Note also that 
$$\ell(r_{i_1} r_{i_2} \dots r_{i_l})=i_1+i_2+\dots + i_l$$
if $1 \leq i_1 < i_2 < \dots < i_l \leq n$. 
Therefore, in the subset $X_n^{(l)}$, there is a unique element of minimal 
length which will be denoted by $a_l=r_1 r_2 \dots r_l$ (note that $a_0=1$). 
It is also clear that $a_l$ has minimal length among all the elements
of $t$-length equal to $l$ (in $W_n$).

Let us make some further notations. We denote by:
\begin{itemize}
\item $\Sym_l$ the parabolic subgroup of $W_n$ generated by $\Sigma_l$
and by $\sigma_l$ its longest element.
\item $W_l$  the parabolic subgroup of $W_n$ generated by $S_l=\{ t\}\cup
\Sigma_l$ and by $w_l$ its longest element.
\item $\Sym_{l,n{-}l}$ the parabolic subgroup of $W_n$ generated by 
$\Sigma_{l,n{-}l}=S_n\setminus \{ t,s_l\}$ and by $\sigma_{l,n{-}l}$ 
its longest element.
\item $W_{l,n{-}l}$  the parabolic subgroup of $W_n$ generated by 
$S_{l,n{-}l}=S_n\setminus \{ s_l\}$ and by $w_{l,n{-}l}$ 
its longest element.
\item $Y_{l,n{-}l}$ the set of distinguished left coset representatives
of $\Sym_{l,n{-}l}$ in $\Sym_n$.
\end{itemize}
One can immediately check that 
$$a_l=w_{l,n-l} \s_{l,n-l}=\s_{l,n-l} w_{l,n-l}=w_l \s_l = \s_l w_l.$$ 
In particular, $a_l^2=1$ and conjugacy by $a_l$ stabilizes $S_l$, $\Sig_l$, 
$S_{l,n-l}$ and $\Sig_{l,n-l}$. In particular, $a_l$ normalizes $W_l$, 
$\SG_l$, $W_{l,n-l}$ and $\SG_{l,n-l}$. Note also that  
$$|X_n^{(l)}|=\matrice{n \\ l}=\frac{n!}{l! (n-l)!}.$$

One can notice that $|X_n^{(l)}|=|Y_{l,n-l}|$ and that 
$|X_n^{(l)}.\SG_{l,n-l}|=n!=|\SG_n|$. This is not a coincidence, as it 
will be shown in this subsection (see Proposition \ref{bijections xnl}). 
We will need the following elementary lemma (the proof is left to 
the reader).

\begin{lem}\label{egalite ri}
If $1 \leq i \leq n$ then
$$r_i^{-1}(j)=\left\{\begin{array}{ll}
                     j+1 & {\mathit{if}}~1 \leq j \leq i-1 \\
                     -1  & {\mathit{if}}~j=i \\
                     j   & {\mathit{if}}~i+1 \leq j \leq n
		     \end{array}\right.$$
Therefore, if $1 \leq j < i-1$, then $r_i^{-1} s_j r_i =s_{j+1}$.
\end{lem}

\begin{cor}\label{conjugaison xnl} 
Let $a \in X_n^{(l)}$ with $1 \leq l \leq n-1$ and let 
$s \in \Sig_n$. Then $a^{-1} s a \not= s_l$.
\end{cor}

\begin{proof} Let us write $a=r_{i_1} \dots r_{i_l}$ with 
$1 \leq i_1 < \dots < i_l \leq n$. 
We want to prove by induction on $l \ge 1$ that 
$a s_l a^{-1} \not\in \Sig_n$.

Let us assume first that $l=1$. We write $i=i_1$. 
Then $(r_i s_1 r_i^{-1})(i)=(r_i s_1)(-1)=r_i(-2)=-r_i(2) < 0$. 
So $r_i s_1 r_i^{-1} \not\in \Sig_n$. 

Let us assume now that $l \ge 2$ and that the 
result holds for $l' < l$. Two cases may occur~:

\medskip

\noindent{\it First case. Assume that $i_l \geq l+1$.} Therefore, 
by Lemma \ref{egalite ri}, we have $r_{i_l} s_l r_{i_l}^{-1}=s_{l{-}1}$ 
and the result follows by induction.

\medskip

\noindent{\it Second case. Assume that $i_l \leq l$.} In this case, 
we have $i_l=l$ and $a=a_l$. But,  
$$a_l s_l a_l^{-1}(l+1)=a_ls_l(l+1)=a_l(l)=-1.$$
So $a_ls_l a_l^{-1} \not\in \Sig_n$.
\end{proof}

\begin{prop}\label{bijections xnl}
Let $l$ be a non-negative integer such that $0 \leq l \leq n$. 
Then~:

\begin{itemize}
\item[(a)] The map $Y_{l,n-l} \to X_n^{(l)}$, $w \mapsto wa_l$ is a bijection.

\item[(b)] The map $\SG_n \to X_n^{(l)}.\SG_{l,n-l}$, $w \mapsto wa_l$ 
is a bijection.
\end{itemize}
\end{prop}

\begin{proof} 
(a) Since $|Y_{l,n-l}|=|X_n^{(l)}|$, it is sufficient to 
prove that $w a_l \in X_n^{(l)}$ 
for every $w \in Y_{l,n-l}$. So let $w \in Y_{l,n-l}$ and
let $s \in \{s_1,s_2,\dots,s_{n-1}\}$. We need to prove that 
$\ell(wa_ls) = \ell(wa_l)+1$. If $a_l s \in (X_n^{(l)})^{-1}$, then 
$\ell(a_l s) = \ell(a_l)+1$ (because $a_l$ is the element of $X_n^{(l)}$ of minimal 
length) and 
$$\ell(wa_ls)=\ell(w)+\ell(a_ls) = \ell(w)+\ell(a_l)+1 = \ell(wa_l)+1.$$ 
If $a_l s \not\in (X_n^{(l)})^{-1}$, then, by Deodhar's Lemma, we have 
$a_l s= s' a_l$ for some $s' \in \{s_1,s_2,\dots,s_{n-1}\}$. But, by Corollary 
\ref{conjugaison xnl}, we have $s' \not= s_l$. So $\ell(ws') = \ell(w)+1$. 
Therefore, 
$$\ell(wa_l s)=\ell(ws'a_l)=\ell(ws')+\ell(a_l) =\ell(w)+\ell(a_l)+1=
\ell(wa_l)+1.$$

(b) Since 
$|X_n^{(l)}.\SG_{l,n-l}|=|\SG_n|$, we only need to show 
that $wa_l \in X_n^{(l)}.\SG_{l,n-l}$ for every $w \in \SG_n$. 
So let $w \in \SG_n$. We write $w=y w'$ with $y \in Y_{l,n-l}$ and 
$w' \in \SG_{l,n-l}$. Then $wa_l=ya_l.a_l^{-1}w'a_l$. But, 
by (a), we have that $ya_l \in X_n^{(l)}$. Moreover, $a_l$ normalizes 
$\SG_{l,n-l}$, so $a_l^{-1}w'a_l \in \SG_{l,n-l}$. Therefore, 
$wa_l \in X_n^{(l)}.\SG_{l,n-l}$.
\end{proof}

The next result concerns the Bruhat ordering.

\begin{prop}\label{ordre xnl}
Let $x$ and $y$ be two elements of $W_n$ such that $x \leq y$. 
We assume that $\ell_t(x)=\ell_t(y)=l$ and that $y \in X_n^{(l)}. \SG_{l,n-l}$. 
Then $x \in X_n^{(l)}. \SG_{l,n-l}$.
\end{prop}

\begin{proof} 
By Proposition \ref{bijections xnl} (b), there exists $w \in \SG_n$ 
such that $y=wa_l$. So there exists a reduced expression
$x=w'a'$ where $w' \leq w$ and $a' \leq a_l$. 
But $\ell_t(x)=l$ so $\ell_t(a')=l$. Therefore $a'=a_l$. 
Moreover $w' \in \SG_n$. So $x=w' a_l \in X_n^{(l)}. \SG_{l,n-l}$ 
by Proposition \ref{bijections xnl} (b).
\end{proof}
\end{what}

\begin{rem} \label{tableaux Sn}
(a) Let $\pi_1, \pi_2\in\Sym_n$. A simple computation together with
\ref{caracRS} show that we have  \[B(\pi_1)=B(\pi_2)\Longleftrightarrow
B(\sigma_n\pi_1)=B(\sigma_n\pi_2).\]

(b) Let $\sigma =\pi\rho\in\Sym_{l,n{-}l}$, with $\pi\in\Sym_l$
and $\rho\in \SG_{[l+1,n]}$ (where $\Sym_{[l+1,n]}$ is
the parabolic subgroup of $\Sym_n$ generated by $s_{l{+}1},\ldots ,
s_{n{-}1}$). Let $(A(\pi), B(\pi))$ and $(A(\rho), B(\rho))$  
the pairs of standard tableaux associated by the (classical) 
Robinson-Schensted correspondence to $\pi$ and $\rho$ respectively.
Then the pair of tableaux associated to $\sigma\in\Sym_n$ is 
given by $(A(\pi)\cdot A(\rho), B(\pi)\cdot B(\rho))$. (For the 
definition of the product of tableaux see \cite[1.2]{Fulton}.)
\end{rem}

\begin{what}\label{decomposition abc}
{\bf A decomposition of elements of $W_n$.} 
If $w \in W_n$ and if $l=\ell_t(w)$, then, by Proposition 
\ref{bijections xnl} (a), there exist uniquely determined 
elements $a_w$, $b_w \in Y_{l,n-l}$, and $\s_w \in \SG_{l,n-l}$ 
such that $$w=a_w a_l \s_w b_w^{-1}.$$
Moreover, we have
$$\ell(w)=\ell(a_w)+\ell(a_l)+\ell(\s_w)+ \ell(b_w).$$
If $w$ and $w'$ are two elements of $W_n$, we write $w \cedricl w'$ 
if 
$$\ell_t(w)=\ell_t(w'), \qquad\s_w \lacril \s_{w'}\qquad{\mathrm{and}}\qquad
b_w=b_{w'}.$$
It is obvious that $\cedricl$ is an equivalence relation. We write 
$\s_l\s_w=\s_w^\prime\s_w^{\prime\prime}$ with $\s_w^\prime \in \SG_l$ 
and $\s_w^{\prime\prime} \in \SG_{[l+1,n]}$. Then 
$$w=a_w w_l \s_w^\prime \s_w^{\prime\prime} b_w^{-1}.$$
Moreover, $w_l$, $\s_w^\prime$ and $\s_w^{\prime\prime}$ commute to each other. 
By Remark \ref{tableaux Sn} we have $w \cedricl w'$ if and only if 
$$\ell_t(w)=\ell_t(w'), \quad\s_w^\prime \lacril \s_{w'}^\prime, \quad 
\s_w^{\prime\prime} \lacril \s_{w'}^{\prime\prime}\quad{\mathrm{and}}\quad
b_w=b_{w'}.$$

\begin{what} \label{shapelm} {\bf Bitableaux and decomposition}.
Note that by the (usual) Robinson-Schensted correspondence, the element
$\sigma ' _w\sigma '' _w$ is associated to a pair of $n$-standard 
bitableaux of shape $(\lambda ,\mu )$ with $|\l |=l$ and $|\m |= n-l$; 
also note that (by Proposition 
\ref{caracRS}) the elements $a_w$ and $b_w$ do not affect the shape of 
the bitableaux of $w_l\sigma ' _w\sigma '' _w$. This implies that the 
shape of the $n$-standard bitableaux associated to $w$ by the generalized 
Robinson-Schensted correspondence is exactly $(\mu ,\lambda )$.
As for the filling of the $n$-standard bitableaux associated to $w$:
$A_n^{+}$ is obtained by the action of $a_w$ on the $A$-tableau of
$\sigma ''$, $A_n^{-}$ is obtained by the action of $a_w$ on the
$A$-tableau of $\sigma '$, while $B_n^{+}$ is obtained by the action of 
$b_w$ on the $B$-tableau of $\sigma ''$ and finally $B_n^{-}$ is obtained 
by the action of $b_w$ on the $B$-tableau of $\sigma '$.  
\end{what}

The above remark has as direct consequence the next proposition.
\begin{prop}\label{equivalent}
Let $w$ and $w'$ be two elements of $W_n$. Then the following are equivalent~:

$(1)$ $w \lacril w'$~;

$(2)$ $w \cedricl w'$.
\end{prop}

\begin{proof} Follows immediately from Proposition \ref{caracRS} and
Remark \ref{shapelm}.
\end{proof}

\begin{rem} \label{involtns} Note that Proposition \ref{equivalent} implies 
that every equivalence class for $\cedricl$ contains a unique involution.
\end{rem}

We end this subsection with a result about the above 
decomposition of elements of $W_n$ and the Bruhat order. 

\begin{prop}\label{abc bruhat}
Let $x$ and $y$ be two elements of $W_n$ such that $\ell_t(x)=
\ell_t(y)=l$ and $x \leq y$. Then~:
\begin{itemize}
\item[(a)] $a_x \leq a_y$ and $b_x \leq b_y$.
\item[(b)] If $y \in X_n^{(l)}. \SG_{l,n-l}$, then $x \in X_n^{(l)}
.\SG_{l,n-l}$.
\item[(c)] If $b_x=b_y=b$, then $xb \leq yb$ and the map 
$[xb,yb] \to [x,y]$, $z \mapsto zb^{-1}$ is an isomorphism between
the two intervals.
\end{itemize}
\end{prop}

\begin{proof} 

(a) By \cite[Lemma 9.10 (f)]{Lusztig02}, we have $a_x a_l \leq a_y a_l$. 
But, $a_l$ is minimal in $\SG_n a_l$. Therefore, again by 
\cite[Lemma 9.10 (f)]{Lusztig02}, we get that $a_x \leq a_y$. 

On the other hand, we have $x^{-1} \leq y^{-1}$. So, by the previous 
result, we have $a_{x^{-1}} \leq a_{y^{-1}}$. But Proposition 
\ref{bijections xnl} together with a simple computation show that 
$a_{x^{-1}}=b_x$. Hence $b_x \leq b_y$. 

(b) By Proposition \ref{bijections xnl} (b), there exists $w \in \SG_n$
such that $y=wa_l$. So there exists a reduced expression
$x=w'a'$ where $w' \leq w$ and $a' \leq a_l$.
But $\ell_t(x)=l$ so $\ell_t(a')=l$. Therefore $a'=a_l$.
Moreover $w' \in \SG_n$. So $x=w' a_l \in X_n^{(l)}. \SG_{l,n-l}$
by Proposition \ref{bijections xnl} (b).

(c) Let us write $x=\a\b$ (reduced expression) where $\a \leq yb$
and $\b \leq b^{-1}$. But, by (b), we have $\a \in X_n^{(l)}.
\SG_{l,n-l}$. So $\b=b^{-1}$ and $xb=\a \leq yb$. 

Now, let us prove that the map 
$f : [xb,yb] \to [x,y]$, $z \mapsto zb^{-1}$ is an increasing bijection. 
First note that the map is well-defined, injective and increasing. 
We need to show that it is surjective. Let $z \in [x,y]$. By (a), we have 
$b=b_x \leq b_z \leq b_y=b$. So $b_z = b$. In particular, the previous 
result shows that $xb \leq zb \leq yb$. Therefore, $zb \in [xb,yb]$ 
and $f(zb)=z$. 
Let $u,v\in [x,y]$ such that $u\leq v$; then there exist $u',v'\in
[xb,yb]$ such that $u=u'b^{{-}1}$ and $v=v'b^{{-}1}$. Now $u'\leq v'$
for the same reasons for which $xb\leq yb$.
\end{proof}

\end{what}

\section{Kazhdan--Lusztig polynomials in the unequal parameter case} 
\label{sec4}

In this section, we recall the basic constructions from Lusztig 
\cite{Lusztig83}. Let $(W,S)$ be a Coxeter system and $\ell \colon W
\rightarrow {\N}_0$ the corresponding length function. Let $\varphi \colon 
W\rightarrow \Gamma$ be a map into an abelian group $\Gamma$ such that 
$\varphi(w_1w_2)= \varphi(w_1)\varphi(w_2)$ whenever $\ell(w_1w_2)=\ell(w_1)+
\ell(w_2)$. Note that this condition implies that $\varphi$ is determined
by its images on $S$ and that $\varphi (s)=\varphi (s')$ 
whenever $s, s'\in S$ are conjugate in $W$. We set 
\[ \varphi(w)=v_w \qquad \mbox{for any $w\in W$}.\]
Let $\cH$ be the generic Iwahori-Hecke algebra associated with
$(W,S)$ over the ring $A={\Z}[\Gamma]$; then $\cH$ has a 
basis $\{T_w\mid w\in W\}$ such that the multiplication is given by
\[T_s^2=T_1+(v_s-v_s^{-1})T_s \qquad \mbox{for $s\in S$}\]
and $T_{w_1}T_{w_2}=T_{w_1w_2}$ whenever $\ell(w_1w_2)=\ell(w_1)+\ell(w_2)$.
Each $T_s$ ($s \in S$) is invertible in $\cH$; we have
\[ T_s^{-1}=T_s+(v_s^{-1}-v_s)T_1 \qquad \mbox{for all $s \in S$}.\]
(note that this is the basis usually denoted by 
$\{ \tilde{T}_w\mid w\in W\}$, see \cite{Lusztig83}).
The construction of a Kazhdan--Lusztig basis of $\cH$ depends on one
further ingredient, namely, the choice of a total ordering on $\Gamma$ which
is compatible with the group structure. Thus, we assume that we have 
fixed a multiplicatively closed subset $\Gamma_{+} \subset \Gamma$ such 
that we have a disjoint union $\Gamma=\Gamma_{+} \cup \{1\}\cup\Gamma_{-}$, 
where $\Gamma_{-}=\{g^{-1} \mid g \in \Gamma_{+}\}$. Note that this
means, in particular, that $1_{\Gamma} \not\in \Gamma_{+}$. We assume that 
\[v_s \in \Gamma_{+} \qquad \mbox{for all $s \in S$}.\]

\begin{exmp} \label{expkl} Let $\Gamma$ be the infinite cyclic group 
generated by some indeterminate~$v$ over $\C$, and $\Gamma_{+}=
\{v^m \mid m  \geq 1\}$. Let $\{c_s\mid s\in S\}$ be a collection of 
positive integers such that $c_s=c_t$ whenever $s$, $t\in S$ are conjugate 
in $W$.  Then the above requirements are satisfied for the unique map
$\varphi \colon W\rightarrow \Gamma$ such that $\varphi(s)=v^{c_s}$
for all $s\in S$. \finl
\end{exmp}

\begin{rem} \label{rem1} 
One should keep in mind that the choice of 
$\varphi \colon W\rightarrow \Gamma$ in Example~\ref{expkl} is the most 
important one as far as applications to representations of reductive 
algebraic groups are concerned; see \cite{Lusztig84} and \cite{Lusztig90b}.
However, more general choices of $\varphi\colon W\rightarrow \Gamma$ have 
applications to the representation theory of Iwahori--Hecke algebras, via 
the construction of left cell representations. For example, the 
determination of the left cells for a two-parameter algebra of type $F_4$ 
in \cite[Chap.~11]{ourbuch} has lead to a construction of the irreducible 
representations in terms of $W$-graphs in this case. These $W$-graphs in 
turn yield a complete set of irreducible representations for any 
semisimple specialisation of that algebra. 
\end{rem}

\begin{what} \label{klbasis} {\bf The Kazhdan--Lusztig basis.}
Let $a \mapsto \bar{a}$ be the involution of ${\Z}[\Gamma]$ which takes 
$g$ to $g^{-1}$ for any $g \in \Gamma$. We extend it to a map $\cH 
\rightarrow \cH$, $h \mapsto \overline{h}$, by the formula 
\[ \overline{\sum_{w \in W} a_w T_w}=\sum_{w \in W} \bar{a}_w 
T_{w^{-1}}^{-1} \qquad (a_w \in {\Z}[\Gamma]).\] 
Then $h \mapsto \overline{h}$ is in fact a ring involution. In this 
set-up, let $\{C_w\mid w \in W\}$ be the basis of $\cH$ constructed in 
\cite[Prop.~2]{Lusztig83} (formerly denoted by $\{ C_w'\mid w\in W\}$).
 We have  
\begin{equation*}
C_w=\sum_{\atop{y\in W}{y \leq w}} P_{y,w}^* T_y
\qquad \mbox{with $P_{y,w}^* \in A$},\tag{a}
\end{equation*}
where $P_{w,w}^*=1$ and $P_{y,w}^*\in {\Z}[\Gamma_{-}]$ if $y< w$. Here,
${\Z}[\Gamma_-]$ denotes the set of integral linear combinations of 
elements in $\Gamma_-$. 

We have $\overline{P_{y,w}^{*}}-
P_{y,w}^{*}=\displaystyle\sum_{y<x\leq w}R_{y,x}P_{x,w}^{*}$,
where the coefficients $R_{x,w}\in\Z [\Gamma ]$ are defined by
$\overline{T_w}=\displaystyle\sum_{y\leq w}R_{y,w}T_y$.
Note also that 
\begin{equation*}
P_{y^{-1},w^{-1}}^*=P_{y,w}^*\qquad \mbox{for all $y,w\in W$ such that 
$y\leq w$}.  \tag{b}
\end{equation*}
This immediately follows from the discussion in \cite[\S 6]{Lusztig83}.
\end{what}

\begin{what}\label{mpoly} {\bf The $M$-polynomials.} Let $w \in W$ 
and $s \in S$ be such that $sw>w$. As in \cite[\S 3]{Lusztig83}, for each 
$y \in W$ such that $sy<y<w<sw$ we define an element $M_{y,w}^s \in 
{\Z}[\Gamma]$ by the inductive condition
\begin{equation*}
M_{y,w}^s + \sum_{\atop{y < z <w}{sz<z}} P^{*}_{y,z}M_{z,w}^s-v_s 
P^{*}_{y,w} \in {\Z}[\Gamma_{-}] \tag{M1} \end{equation*}
and the symmetry condition 
\begin{equation*} \overline{M}_{y,w}^s=M_{y,w}^s.\tag{M2}
\end{equation*}

With the above definition we have the following multiplication formulas; 
see \cite[Prop.~4]{Lusztig83}. Let $w \in W$ and $s \in S$. Then we have 
\[ \renewcommand{\arraystretch}{2}
C_{s}\, C_{w}= \left\{ \begin{array}{ll}
C_{sw}+\displaystyle\sum_{\atop{z<w}{sz<z}}M_{z,w}^s\, C_z \qquad
\text{if}\quad sw>w\\
(v_s+v_s^{{-}1})\, C_{w} \qquad \text{if} \quad sw<w\, .
\end{array} \right. \]

The proof for the above multiplication rules in \cite{Lusztig83} actually
provides a recursion formula for the computation of $P^{*}_{y,w}$. 
First recall that $P^{*}_{w,w}=1$ for all $w\in W$ and $P^{*}_{y,w}=0$
unless $y\leq w$. Now let $y$ and $w$ be two elements of $W$ 
such that $y<w$ and let $s \in S$ such that $sw < w$. Then:
\begin{alignat}{2}\label{p1p2} 
P^{*}_{y,w} &= v_sP^{*}_{y,sw}+P^{*}_{sy,sw}- 
\sum_{\atop{y \leq z<sw}{sz<z}} P^{*}_{y,z}M_{z,sw}^s &&\qquad 
\mbox{if $sy<y$},\\ P^{*}_{y,w} &= v_s^{-1}P^{*}_{sy,w} &&\qquad 
\mbox{if $sy>y$}.
\end{alignat}

We conclude by an obvious lemma concerning the degree of the $M$-polynomials.

\begin{lem} \label{lemmeM} Let $s\in S$ and $y,w\in W$ be such 
that $sy<y<w<sw$. Then $v_s^{-1}~ M_{y,w}^{s}\in\Z 
[\Gamma _{-}]$.
\end{lem}

\begin{proof} 
If $sy < y < w < sw$, we set $\Mti_{y,w}^s=v_s^{-1}~M_{y,w}^s$. Then condition (M1) 
implies that 
$$\Mti_{y,w}^s + \sum_{\atop{y < z <w}{sz<z}} P^{*}_{y,z}\Mti_{z,w}^s-
P^{*}_{y,w} \in {\Z}[\Gamma_{-}].$$
Since $P_{y,w}^* \in \ZM[\G_-]$, the result follows immediately by induction 
on $\ell(w)-\ell(y)$.
\end{proof}
\end{what}

\begin{what} {\bf The longest element.} \label{longest} Assume that 
$W$ is finite and let $w_0\in W$ be the unique element of maximal length.
Then, for $y \leq w$ in $W$, we have the following relation 
\begin{equation*}
\sum_{\atop{z \in W}{y \leq z \leq w}} (-1)^{\ell(w)-\ell(z)} P_{y,z}^*\,
P_{ww_0,zw_0}^*=\left\{\begin{array}{cl} 1 & \quad \mbox{if $y=w$},\\0 & 
\quad \mbox{if $y<w$}.\end{array}\right. \tag{a}
\end{equation*}
Furthermore, if $y,w\in W$ and $s \in S$ are such that $sy<y<w<sw$, then
\begin{equation*}
M_{ww_0,yw_0}^s=-(-1)^{\ell(w)-\ell(y)}M_{y,w}^s.\tag{b}
\end{equation*}
In the equal parameter case, these relations are already contained in 
Kazhdan and Lusztig \cite{KaLu}; for the general case, see \cite[\S 2]{Ge1}.
Passing to the polynomials $P_{y,w}=v_y^{-1}v_w P^*_{y,w}$, 
we obtain
\begin{equation*}
\sum_{\atop{z \in W}{y \leq z \leq w}} (-1)^{\ell(w)-\ell(z)} P_{y,z}\, 
P_{ww_0,zw_0}=\left\{\begin{array}{cl} 1 & \quad \mbox{if $y=w$},\\0 & 
\quad \mbox{if $y<w$}.\end{array}\right. \tag{c}
\end{equation*}
The last relation will be helpful in the computation of certain 
Kazhdan--Lusztig polynomials in type $B_n$.
\end{what}

\begin{what} {\bf Left cells.} \label{cells} We recall the definition
of left cells from \cite[\S 8.1]{Lusztig02}. Let $\leq_L$ be the preorder
relation on $W$ generated by the relation:
\begin{equation*}
\left\{\begin{array}{ll} 
\mbox{$y \leq_{L,s} w$ if there exists some $s \in S$
such that $C_y$ appears with}\\\mbox{non-zero coefficient
in $C_sC_w$ (expressed in the $C_w$-basis)}.
\end{array}\right.\tag{L}
\end{equation*}
The equivalence relation associated with $\leq_L$ will be denoted by $\sim_L$
and the corresponding equivalence classes are called the {\em left cells} of
$W$. Lusztig \cite[\S 8.3]{Lusztig02} has associated to each left cell $\fC$ 
a representation of~$\HC$. 
\end{what}

\begin{rem} When $(W,S)$ is the Coxeter group of type $A_{n{-}1}$,
Kazhdan and Lusztig \cite[\S 5]{KaLu} proved that two elements
$y,w\in W$ are in the same left cell if and only if they are in the 
same Robinson-Schested cell (i.e.~$B(y)=B(w)$). 
\end{rem}

\begin{exmp} \label{b3lusz} Assume that $(W,S)=(W_n,S_n)$ as defined 
in \S\ref{setup}. Let $v$ be an indeterminate and $\Gamma:=
\{v^n \mid n \in \Z\}$. Let $c$, $d \geq 1$. As in Example~\ref{expkl}, 
let $\varphi \colon W\rightarrow \Gamma$ be defined by 
\[ \varphi(t)=v^c \mbox{ and } \varphi(s_i) =v^d 
\mbox{ for $1\leq i \leq n-1$}.\]
Let $\Gamma_+=\{v^n \mid n>0\}$. Assume that $(c,d)\in \{(1,2),(3,2)\}$.
Then Lusztig \cite[Theorem~11]{Lusztig83} has shown that the corresponding 
left cell representations are all irreducible. As an example, we give the
distribution of the elements of $W_3$ into left cells in Table~\ref{tab2}.
(The computation were done using the CHEVIE system \cite{chevie2}.)
We only point out here that the generalized Robinson--Schensted
cells in Table~\ref{tab1} appear to be completely unrelated to the
left cells in Table~\ref{tab2}. Our main result, which will be proven
in Section~\ref{sec7} will show that the Robinson--Schensted cells are
related to the case where $d=1$ and $c$ is large enough~; we shall call
it {\em the asymptotic case}.
\end{exmp}

\begin{what} {\bf Left cells and parabolic subgroups.} \label{geck result}
If $J$ is a subset of $S$, we denote by $W_J$ the parabolic subgroup
of $W$ generated by $J$ and by $X_J$ the set of
distinguished left coset representatives of $W_J$ in $W$.

The next result, due to Geck \cite{geck xi}, will be used in the sequel.

\begin{thm} \label{parabolicG}
{\em [Geck \protect{\cite[Th.~1, proof of Th.~1]{geck xi}}]}
\begin{itemize}
\item[(a)] Let $\mathfrak{C}$ be a left cell of $W_J$. Then
$X_J\cdot\mathfrak{C}$ is a union of left cells of $W$.
\item[(b)] Let $z,u\in X_J$ and $x,y\in W_J$. Then we have
\begin{align*}
& zx\leq_L uy \;\Longrightarrow\; x\leq_L y \quad\text{in}\quad
W_J\quad\text{and}\\ 
& zx\sim_L uy \;\Longrightarrow\; x\sim_L y \quad\text{in}\quad
W_J\, .
\end{align*}
\end{itemize}
\end{thm}
\end{what}

\begin{table}[htbp]\caption{Left cells in type $B_3$ with unequal 
parameters} \label{tab2}
\begin{center}
$\renewcommand{\arraystretch}{1.2} 
\begin{array}{c} \hline \mbox{Parameters $q,q^2,q^2$}\\\hline
\{1\}, \{t\}, \{s_1ts_1s_2s_1ts_1s_2\}, \{ts_1ts_1s_2s_1ts_1s_2\},
\{s_1ts_1,s_2s_1ts_1\},\\ \{ts_1ts_1,ts_2s_1ts_1\}, 
\{s_1ts_1s_2,s_2s_1ts_1s_2\},
\{ts_1ts_1s_2,ts_2s_1ts_1s_2\},\\
\{s_1,s_2s_1,ts_1\}, \{s_1t,s_2s_1t,ts_1t\},
\{s_2,s_1s_2,ts_1s_2\}, \{ts_2,s_1ts_2,ts_1ts_2\},\\
\{ts_2s_1,s_1ts_2s_1,ts_1ts_2s_1\}, \{ts_2s_1t,s_1ts_2s_1t,ts_1ts_2s_1t\},\\
\{s_1s_2s_1,ts_1s_2s_1,s_1ts_1s_2s_1\},
\{s_1s_2s_1t,ts_1s_2s_1t,s_1ts_1s_2s_1t\},\\
\{s_1s_2s_1ts_1,ts_1s_2s_1ts_1,s_1ts_1s_2s_1ts_1\},\\
\{s_1s_2s_1ts_1s_2,ts_1ts_1s_2s_1,ts_1s_2s_1ts_1s_2\},\\
\{s_1ts_2s_1ts_1,ts_1ts_2s_1ts_1,ts_1ts_1s_2s_1ts_1\},\\
\{s_1ts_2s_1ts_1s_2,ts_1ts_1s_2s_1t,ts_1ts_2s_1ts_1s_2\}\\[0.5cm]\hline
\mbox{Parameters $q^3,q^2,q^2$}\\\hline
\{1\}, \{ts_1t\},\{s_1s_2s_1ts_1s_2\},\{ts_1ts_1s_2s_1ts_1s_2\},
\{s_1,s_2s_1\}, \{s_2,s_1s_2\}, \\
\{t,s_1t,s_2s_1t\},\{ts_2,s_1ts_2,ts_1ts_2\},\{ts_1,s_1ts_1,s_2s_1ts_1\},\\
\{ts_2s_1,s_1ts_2s_1,ts_1ts_2s_1\},
\{ts_1s_2,s_1ts_1s_2,s_2s_1ts_1s_2\},\\
\{ts_2s_1t,s_1ts_2s_1t,ts_1ts_2s_1t\},
\{ts_1ts_2s_1ts_1,ts_1ts_1s_2s_1ts_1\},\\
\{s_1s_2s_1,ts_1s_2s_1,s_1ts_1s_2s_1\},
\{ts_1ts_1s_2s_1t,ts_1ts_2s_1ts_1s_2\},\\
\{ts_1ts_1,ts_2s_1ts_1,s_1ts_2s_1ts_1\},
\{s_1s_2s_1t,ts_1s_2s_1t,s_1ts_1s_2s_1t\},\\
\{ts_1ts_1s_2,ts_2s_1ts_1s_2,s_1ts_2s_1ts_1s_2\},\\
\{s_1s_2s_1ts_1,ts_1s_2s_1ts_1,s_1ts_1s_2s_1ts_1\},\\
\{ts_1ts_1s_2s_1,ts_1s_2s_1ts_1s_2,s_1ts_1s_2s_1ts_1s_2\}\\\hline
\end{array}$
\end{center}
\end{table}
\section{Type $B_n$ in the asymptotic case} \label{sec5}

Now consider the group $W_n$, with generators $S_n=\{t,s_1,
\ldots,s_{n-1}\}$ as in Section~\ref{setup}. 
Let $A={\Z}[V,V^{-1},v,v^{-1}]$ where
$V$ and $v$ are independent indeterminates and  
$\Gamma=\{V^iv^j \mid i,j \in \Z\}$ (which is an abelian group
under multiplication). We define $\varphi \colon W_n\rightarrow 
\Gamma$ by 
\[ v_t=V\quad\text{and}\quad v_{s_i}= v\quad\text{for}
\quad 1\leq i\leq n-1\, .\]
In particular, we have 
$\varphi(w)=v_w=V^{\ell_t(w)} v^{\ell_s(w)}$ for $w\in W_n$.
Let $\cH_n$ be the corresponding generic Iwahori--Hecke 
algebra over $A$, with quadratic relations 
\begin{align*}
T_t^2 &=T_1+(V-V^{-1})T_t,\\
T_{s_i}^2 &=T_1+(v-v^{-1})T_{s_i} \quad \mbox{for $1\leq i\leq n-1$}.
\end{align*}

\bigskip 
\bigskip

\begin{quotation}
{\bf Hypothesis :} 
{\em We fix a lexicographic ordering of $\Gamma$ where }
\begin{equation*}
\Gamma_+=\{V^i v^j \mid i>0, \mbox{ any $j$}\}\cup \{v^j\mid
j>0\}.\tag{H}
\end{equation*}
\end{quotation}

\bigskip
\bigskip

We have a corresponding Kazhdan--Lusztig basis $\{C_w\}$ in $\cH_n$,
corresponding Kazhdan--Lusztig polynomials $P_{y,w}$, and corresponding
left cells in $W_n$. Note that all these depend on the choice of
$\Gamma_+$. We begin with the following remark.

\begin{rem} \label{independent} In the setting of Example~\ref{b3lusz},
assume that $d=1$. Then there exists a constant $c_0> 1$ such that, for
all $c\geq c_0$, the corresponding left cells of $W_n$ are precisely the 
left cells with respect to the ordering~($*$).

(This follows from the fact that we only have a finite list of 
Kazhdan--Lusztig polynomials and $M$-polynomials for our given group $W_n$.
All these polynomials are two-variable Laurent polynomials in $v$ 
and $V$. Hence it is clear that if we specialize $V$ to a 
sufficiently large power of $v$, then the specialisations of all those 
polynomials will remain in $\ZM[\G_-]$.)

For example, taking $c_0 =n^{2}{-}1$ not only leads to the same left
cells, but also the Kazhdan-Lusztig basis $\{ C_w \}$ of $\cH_n$
specializes (via $V\mapsto v^c$) to the Kazhdan-Lusztig basis in the 
setting of Example~\ref{b3lusz}. Furthermore the multiplication
polynomials $M_{y,w}^s$ in $\cH_n$ specialize to the multiplication
polynomials in the setting of \ref{b3lusz}.

It is likely that there is a better bound, but this one above ($c_0
=n^2{-}1$) will suffice for our purposes.\finl
\end{rem}

\begin{what} {\bf Some properties of the polynomials $M_{y,w}^{s}$.}
We will now establish some basic properties of the polynomials 
$P_{y,w}$ and of the $M$-polynomials. We set $Q=V^2$ and $q=v^2$.

\begin{thm} \label{mythm} With the above choice of $\Gamma_+\subset
\Gamma$, the following hold.
\begin{itemize}
\item[(a)] For all $y,w\in W_n$ with $y\leq w$, we have $P_{y,w}
\in {\Z}[q]$. The constant term of $P_{y,w}$ is $1$.
 
\item[(b)] Let $y,w\in W_n$ and $1\leq i \leq n-1$ be such that
$s_iy<y<w<s_iw$. Then we have $M_{y,w}^{s_i}=0$ unless $\ell_t(y)=\ell_t(w)$.
Furthermore, if $\ell_t(y)=\ell_t(w)$, then $P_{y,w}\in {\Z}[q]$ has degree 
at most $(\ell(w)-\ell(y)-1)/2$ and $M_{y,w}^{s_i}\in \Z$ is the coefficient of 
$q^{(\ell(w)-\ell(y)-1)/2}$. 
\end{itemize}
\end{thm}

\begin{proof} Let us make a preliminary remark. By combining the symmetry
condition (M2) in the definition \ref{mpoly} of $M$-polynomials with 
Lemma \ref{lemmeM} and with our choice for $\Gamma_{+}$ we obtain 
that $M_{y,w}^{s_i}\in\Z$, for all  $y,w\in W_n$ and $1\leq i \leq n-1$ 
such that $s_iy<y<w<s_iw$. 

(a) We will proceed by induction.
For this purpose, it will be convenient to consider all groups $W_n$
at the same time. Note that we have standard embeddings $W_0\subset W_1 
\subset W_2 \subset W_3 \subset \cdots$, where $W_0=\{1\}$. Furthermore, 
if $y,w\in W_n$ lie in $W_m$ for some $m<n$, then $P_{y,w}$ computed 
with respect to $W_m$ is the same as $P_{y,w}$ computed with respect 
to $W_n$; a similar result also holds for the $M$-polynomials. (This 
immediately follows, for example, from the recursion formulas  
for the Kazhdan-Lusztig polynomials in \cite[11.1]{ourbuch}.) 
Now let $W:=\bigcup_{n\geq 0}W_n$. For a pair 
$(y,w)$ of elements in $W$ such that $y\leq w$, we set 
\[ \lambda(y,w)=(n(w),\ell(w)-\ell(y), \ell_t(w),\ell_s(w))\]
where $n(w):=\min\{n \geq 0 \mid w \in W_n\}$. Let $\infspe$ be the 
usual lexicographical ordering of these quadruples. We write
$\lambda(y',w')\inferieur \lambda(y,w)$ if $\lambda(y',w')\infspe
\lambda(y,w)$ and $\lambda(y',w')\neq \lambda(y,w)$.

Now let $(y,w)$ be a pair of elements in $W$ such that $y\leq w$. If
$\lambda(y,w)=(0,0,0,0)$, then $y=w=1$ and there is nothing to prove. 
Now assume that $\lambda(y,w)\neq (0,0,0,0)$ and that the assertions 
hold for all pairs $(y',w')$ of elements in $W$ such that 
$\lambda(y',w')\inferieur \lambda(y,w)$. Let $n=n(w)$. First we show that
$P_{y,w}\in {\Z}[q]$. We distinguish several cases.

{\bf Case 1}. Suppose there exists some $1\leq i \leq n-1$ such that
$s_iw<w$. Then we can apply the recursion formulas  
and see that $P_{y,w}\in {\Z}[q]$ by induction. Note that each term in 
that formula involves a Kazhdan--Lusztig polynomial or an $M$-polynomial 
associated with a pair $(y',w')$ such that $n(w')\leq n(w)$, 
$\ell(w')-\ell(y')\leq \ell(w)-\ell(y)$, $\ell_t(w')=\ell(w)$ and 
$\ell_s(w') \leq \ell_s(w)$,
where at least one inequality is strict. Thus, we have $\lambda(y',w')
\inferieur \lambda(y,w)$. 

{\bf Case 2}. Suppose there exists some $1\leq i \leq n-1$ such that
$ws_i<w$. Then we use the fact that $P_{y,w}=P_{y',w'}$ where $y'=y^{-1}$ 
and $w'=w^{-1}$; see (\ref{klbasis})(c). We have $s_iw'<w'$ and so we can 
apply the argument in Case~1 to conclude that $P_{y,w}\in {\Z}[q]$.

{\bf Case 3}. Suppose that $s_iw>w$ and $ws_i>w$ for all $1\leq i\leq n-1$.
Then the only remaining possibility is that $tw<w$ and $wt<w$. This means
that $w$ is a ``bigrassmannian'' in the sense of Geck and Kim \cite{GeKi}. 
So we have 
\[ w=r_1r_2\cdots r_k \qquad 
\mbox{where $r_i$ is defined as in \S\ref{preliminaires abc}};\]
see \cite[Prop.~4.2]{GeKi}. In particular, since we assumed that $n=n(w)$, 
we must have $k=n$. If $ty>y$, we have $P_{y,w}=P_{ty,w}$ by 
the recursion formulas for $P_{y,w}$. Hence we can apply induction since 
\[\lambda(ty,w) =(n(w),\ell(w)-\ell(y)-1,\ell_t(w),\ell_s(w))\prec \lambda(y,w).\]
So we have $P_{y,w}\in {\Z}[q]$ in this case. Let us now assume that $ty<y$. 
Using the relation in (\ref{longest})(c), we can express $P_{y,w}$ as 
follows.
\[P_{y,w}=-(-1)^{\ell(w)-\ell(y)}P_{ww_n,yw_n}-\sum_{\atop{z \in W}{y<z<w}} 
(-1)^{\ell(w)-\ell(z)} P_{y,z}\, P_{ww_n,zw_n};\]
here, $w_n$ is again the longest element in $W_n$.
Let $(y',w')$ be a pair of elements in $W_n$ occuring in the above sum,
that is, we have $(y',w')=(y,z)$ or $(y',w')=(ww_n,zw_n)$ where
$y<z<w$. Then we have $n(w')\leq n(w)$ and $\ell(w')-\ell(y')<
\ell(w)-\ell(y)$ and so $\lambda(y',w')\inferieur \lambda(y,w)$. 
Thus, by induction, all terms
in the above lie in ${\Z}[q]$. So it remains to check if we can also
apply induction to the term $P_{ww_n,yw_n}\in {\Z}[q]$. This is indeed
so. For, we certainly have $n(yw_n)\leq n(w)$ and $\ell(yw_n)-\ell(ww_n)=
\ell(w)-\ell(y)$. Next note that $\ell_t(w)=n$. On the other hand, 
we have $\ell_t(yw_n)=n-\ell_t(y)<n$ since $ty<y$ by 
assumption. Thus, we have $\lambda(ww_n,yw_n)\inferieur 
\lambda(y,w)$ and so we can apply induction to $P_{ww_n,yw_n}$.

(b) Let us prove the remaining assertions in (b).
We have already seen that $P_{y,w}\in {\Z}[q]$. Thus, if
$\ell_t(y)<\ell_t(w)$, then $vP_{y,w}^*$ is a negative power
of $V$ times  a Laurent polynomial in $v$ and so $vP_{y,w}^*
\in {\Z}[\Gamma_-]$ by the definition of $\Gamma_+$. Hence we have
$M_{y,w}^{s_i}=0$ in this case. Now assume that $\ell_t(y)=\ell_t(w)$.
Since $P_{y,w}\in {\Z}[q]$, we have
$P_{y,w}^*\in {\Z}[q^{-1}]$ if $\ell(w)-\ell(y)$ is even. This shows that
$M_{y,w}^{s_i}=0$ whenever $\ell(w)-\ell(y)$ is even. 

Finally, the fact that $P_{y,w}^*=v^{\ell(y)-\ell(w)}P_{y,w}$ is a 
polynomial in $v^{-1}$ with constant term zero also implies that 
$\deg_q (P_{y,w})\leq (\ell(w)-\ell(y)-1)/2$ and that $M_{y,w}^{s_i}$ 
is the coefficient of $v^{\ell(w)-\ell(y)-1}$. Thus, (b) is proved. 
\end{proof}

\end{what}

\begin{what}\label{mt subsection}{\bf On the polynomials $M_{y,w}^t$.} 
Our aim here is to prove that, if $y,w\in W_n$ are such that
$ty < y < w < tw$ and
$\ell_t(y) < \ell_y(w)$, then $M_{y,w}^t=0$. 
Before proving this, we need the following lemma.

\begin{lem}\label{petit lemme rigolo}
Assume that $n \geq 2$ and 
let $w \in W_n$ satisfying the four following properties~:
\begin{itemize}
\item[(1)] $ws_i > w$ for every $i \in \{1,2,\dots,n-1\}$.

\item[(2)] $s_i w > w$ for every $i \in \{2,\dots,n-1\}$.

\item[(3)] $tw>w$. 

\item[(4)] $w \not\in W_{n-1}$.
\end{itemize}
\noindent Then $w=r_2r_3 \dots r_n$.
\end{lem}

\begin{proof}
We set $l=\ell_t(w)$. 
By (1), we have $w \in X_n^{(l)}$. So there exist 
$1 \leq i_1 < i_2 < \dots < i_l \leq n$ such that 
$w=r_{i_1} r_{i_2} \dots r_{i_l}$. By (3), we have that $i_1 \geq 2$. 
By (2) we have that $i_j-i_{j-1} = 1$ for every $j \in \{1,2,\dots,l\}$. 
By (4), we have that $i_l=n$. So $l=n-1$ and $w=r_2r_3 \dots r_n$.
\end{proof}

We are now ready to prove the following theorem.

\begin{thm}\label{theo mt}
Let $y$ and $w$ be two elements of $W_n$ such that $ty < y < w < tw$. 
If $M_{y,w}^t \not= 0$, then $\ell_t(y)=\ell_t(w)$.
\end{thm}

\begin{proof}
Assume that the theorem does not hold. Let $n$ be minimal such 
that there exists $y$ and $w$ in $W_n$ such that $ty < y < w < tw$,  
$M_{y,w}^t \not= 0$ and $\ell_t(y)<\ell_t(w)$. 
It is obvious that $n \geq 2$. 
We choose such a pair $(y,w)$ in such a way that $\ell(w)$ 
is minimal. 

Assume first that there exists $i \in \{1,2,\dots,n-1\}$ such that 
$w s_i < w$. We set $w'=ws_i$. By Theorem \ref{mythm} (b), we have
$$C_w = C_{w'} C_{s_i} + \sum_{x < w'~\text{and}~\ell_t(x)=\ell_t(w)} \a_x C_x$$
for some $\a_x \in \ZM[\G]$. Therefore, 
$$C_t C_w = C_t C_{w'} C_{s_i} + 
\sum_{x < w'~\text{and}~\ell_t(x)=\ell_t(w)} \a_x C_t C_x.$$
But, by minimality of $w$, $C_t C_{w'}$ (respectively $C_t C_x$) 
has a non-zero coordinate on $C_z$ only if 
$\ell_t(z)=\ell_t(w')=\ell_t(w)$ (respectively 
$\ell_t(z)=\ell_t(x)=\ell_t(w)$) or if $z=tw'$ (or $z=tx$).
So, by Theorem \ref{mythm} (b), 
$C_t C_{w'} C_{s_i}$ has a non-zero coordinate on $C_z$ 
only if $\ell_t(z)=\ell_t(w')=\ell_t(w)$. Therefore, $M_{y,w}^t \not=0$ 
implies that $\ell_t(y)\geq\ell_t(w)$, which is contrary to our hypothesis.
So we have~:

\medskip

\qquad (1) {\it $ws_i > w$ for every $i \in \{1,2,\dots,n-1\}$.}

\medskip

By a similar argument, and using the fact that $C_t$ and $C_{s_i}$ commute 
if $i \geq 2$, we have~:

\medskip

\qquad (2) {\it $s_iw > w$ for every $i \in \{2,\dots,n-1\}$.}

\medskip

By hypothesis, we have

\medskip

\qquad (3) {\it $tw > w$.}

\medskip

By the minimality of $n$, we have

\medskip

\qquad (4) {\it $w \not \in W_{n-1}$.}

\medskip

Therefore, by (1), (2), (3) and (4) and by Lemma \ref{petit lemme rigolo}, 
we get that $w=r_2r_3\dots r_n=s_1s_2 \dots s_{n-1} r_1 r_2 \dots r_{n-1}$. 
Now, let us write $y=zy'$ with $y' \in W_{n-1}=W_{S_{n-1}}$ and $z \in X_{S_{n-1}}$. 
Note that $s_1s_2 \dots s_{n-1} \in X_{S_{n-1}}$. So, 
by \ref{parabolicG}, we have $y' \leq_L w$. 
Therefore, by the minimality of $n$, we get that 
$\ell_t(y) \geq \ell_t(y') \geq \ell_t(w') =\ell_t(w)$. So $\ell_t(y)=\ell_t(w)$.
\end{proof} 

As a consequence of Theorem \ref{mythm} (b) and of Theorem \ref{theo mt}
we obtain the following statement.

\begin{cor}\label{lt constant}
If $y$ and $w$ are such that $y \leq_L w$, then $\ell_t(y) \geq \ell_t(w)$.
If moreover $y \sim_L w$, then $\ell_t(y)=\ell_t(w)$.
\end{cor}
\end{what}

\section{Main result and consequences} \label{sec7}

We prove in this section the main result of this paper, namely the fact 
that left cells coincide with Robinson-Schensted generalized left cells 
for our choice of parameters (see Theorem \ref{main theo}). 
Note that Corollary \ref{lt constant} 
and Remark \ref{t longueur} are first evidences. 

\begin{what}\label{subsection debut preuve}{\bf Some preliminaries.} 
The next results relate some different Kazhdan-Lusztig polynomials 
by using the decomposition of elements of $W_n$ defined in \S\ref{sec3'}:
$w\in W_n$ is written as $w=a_wa_l\sigma_w b_w^{-1}$ (reduced expression), 
where $l=\ell_t (w),\, a_w,b_w\in Y_{l,n-l},\, \sigma_w\in\Sym_{l,n-l}$.

\begin{prop}\label{left droite}
Let $x$ and $y$ be two elements of $W_n$ such that $x \leq y$, 
$\ell_t(x)=\ell_t(y)$ and $b_x=b_y=b$. Then~:
\begin{itemize}
\item[(a)] $R_{x,y}=R_{xb,yb}$.

\item[(b)] $P_{x,y}^*=P_{x b,y b}^*$.

\item[(c)] $P_{x,y}=P_{x b,y b}$.

\item[(d)] If $s \in S_n$ is such that $sx < x < y < sy$, 
then $M_{x,y}^s=M_{xb,yb}^s$.
\end{itemize}
\end{prop}

\bigskip

\begin{proof}
Note that equalities (a), (b), (c) and (d) for $s\in\Sigma_n$ have 
a meaning because of Proposition \ref{abc bruhat}. It is also easy
to check that $tz<z\Longleftrightarrow tzb<zb$.

Now let us write $b=s_{i_1} s_{i_2} \dots s_{i_r}$ where 
$i_j \in \{1,2,\dots,n-1\}$ and $r=\ell(b)$. Then  
$$R_{a_xa_l\sigma_x,a_ya_l\sigma_y}=R_{a_xa_l\sigma_xs_{i_1},a_ya_l\sigma_y
s_{i_1}}=
R_{a_x\a_l\sigma_xs_{i_1} s_{i_2},a_ya_l\sigma_ys_{i_1}s_{i_2}}=
\dots=R_{x,y}$$
by Lusztig, \cite[Lemma 4.4]{Lusztig02}. This proves (a). 

(b) and (d) follow immediately from (a), from Proposition \ref{abc bruhat} (b) 
and from the fact that $P_{x,y}^*$ and $M_{x,y}^s$ are defined inductively 
using the polynomials $R_{!,?}$. Moreover, (b) clearly implies (c).
\end{proof}

\begin{lem}\label{left b constant}
Let $x$ and $y$ be two elements of $W_n$ such that $\ell_t(x)=\ell_t(y)$ 
and $x \leq_L y$. Then $b_x \leq b_y$. In particular, if $x \sim_L y$, 
then $b_x=b_y$.
\end{lem}

\begin{proof}
By Corollary \ref{lt constant}, we may assume that there exists 
some $s \in S_n$ such that $C_x$ appears with a non-zero 
coefficient in $C_s C_y$. Two cases may occur. If 
$x \leq y$, then $b_x \leq b_y$ by Proposition \ref{abc bruhat} (a). 
Otherwise, we have $x=sy > y$. Since $\ell_t(x)=\ell_t(y)$, this implies 
that $s \not= t$. Therefore $b_x=b_y$ by Proposition 
\ref{bijections xnl} (b).
\end{proof}

\begin{cor}\label{on se fout de b}
Let $x$ and $y$ be two elements of $W_n$ such that $\ell_t(x)=\ell_t(y)$ and 
$b_x=b_y=b$. Then the following are equivalent~:
\begin{itemize}
\item[(1)] $x \leq_L y$.

\item[(2)] $xb \leq_L yb$.
\end{itemize}
\end{cor}

\begin{proof}
This follows immediately from Propositions \ref{abc bruhat} and 
\ref{left droite}.
\end{proof}

\begin{cor}\label{on se fout vraiment de b}
Let $x$ and $y$ be two elements of $W_n$ such that $\ell_t(x)=\ell_t(y)$ and 
$b_x=b_y=b$. Then the following are equivalent~:
\begin{itemize}
\item[(1)] $x \sim_L y$.

\item[(2)] $xb \sim_L yb$.
\end{itemize}
\end{cor}

\end{what}

\begin{what}\label{ca y est on y est}{\bf Left cells in type $B_n$.} 
We are now ready to prove the following result~:

\begin{thm}\label{main theo}
Let $x$ and $y$ be two elements of $W_n$. Then the following are equivalent~:
\begin{itemize}
\item[(1)] $B_n (x)=B_n (y)$, that is, $x$ and $y$ lie in the same
generalized RS-cell.  

\item[(2)] $x \cedricl y$.

\item[(3)] $x \sim_L y$, that is, $x$ and $y$ lie in the same left cell.
\end{itemize}
\end{thm}

\begin{proof}
The equivalence between (1) and (2) has been proved in Theorem 
\ref{equivalent}. 

First, we prove that (3) implies (2). So we assume that $x \sim_L y$. 
Then, by Corollary \ref{lt constant}, we have $\ell_t(x)=\ell_t(y)=l$.  
Moreover, by Lemma \ref{left b constant}, we have $b_x=b_y$. Let $b=b_x=b_y$.  
Then it follows, by Corollary \ref{on se fout vraiment de b}, that
$xb\sim_L yb$. But $xb\, ,yb\in X_n^{(l)}.\Sym_{l,n{-}l}$, and by 
Theorem \ref{parabolicG}(b) we get that $\s_x\sim_L\s_y$. 
Now by \cite[\S 5]{KaLu} 
and \cite{ariki}, this implies that $\s_x \lacril \s_y$. 
Thus $x \cedricl y$ as desired.

We prove that (2) implies (3) by using a counting argument. We have just 
seen that each left cell of $W_n$ is contained in a 
generalized Robinson-Schensted cell, that is 
\[ \#\{\text{left cells}\}\geq \#\{ \text{generalized RS-cells}\}\, .\]
On the other hand, since any irreducible representation of $W_n$ is
realized over $\Q$, it is then well-known that the number of involutions
of $W_n$ equals the number of irreducible direct summands of the
left regular representation of $W_n$. Thus, since the representations
carried by the left cells give a direct sum decomposition of the 
left regular representation, we have
\[\#\{\text{left cells}\}\leq\#\,\text{involutions in $W_n$}=\#\{\text
{generalized RS-cells}\}\] 
\end{proof}
\end{what}

\begin{what}\label{rep}{\bf Characters afforded by left cells representations.} 
Let $K$ be the fraction field of $A$. 
We are now interested in the representation theory of the 
$K$-algebra $\HC_n^K=K \otimes_A \HC_n$. This algebra is split 
semisimple. Let $\LG$ be a left cell in $W_n$ and let $x_0 \in \LG$. 
We denote by $\IC_{\leq_L \LG}$ (resp. $\IC_{<_L \LG}$ the left ideal 
of $\HC_n$ having $(C_x)_{x \leq_L x_0}$ (resp. $(C_x)_{x <_L x_0}$) 
has $A$-basis. We set $\VC_\LG=\IC_{\leq_L \LG}/\IC_{<_L \LG}$. 
This is an $\HC_n$-module which is free over $A$. 

We denote by $A \to \ZM$, $a \mapsto \aha$ the unique morphism 
of $\ZM$-algebras such that $\Vha=\vha=1$. 
Let $\aG$ denote the kernel of this morphism~: it is generated by $v-1$ and $V-1$. 
Then $A/\aG \simeq \ZM$ and the $\ZM$-algebra $\HC_n/\aG \HC_n$ 
may be identified with the group algebra $\ZM W_n$. 
Then $\widehat{\VC}_\LG=\VC_\LG/\aG \VC_\LG$ is a $\ZM W_n$-module 
which is free over $\ZM$. By the computational argument of the proof of 
Theorem \ref{main theo}, we easily get~:

\begin{prop}\label{irreducible}
With the above notation we have~:
\begin{itemize}
\item[(a)] If $\LG$ be a left cell of $W_n$,  
then $K \otimes_A \VC_\LG$ is an irreducible $\HC_n^K$-module 
and $\QM \otimes_\ZM \widehat{\VC}_\LG$ is an irreducible $\QM W_n$-module. 

\item[(b)] Every irreducible $\HC_n^K$-module (resp. $\QM W_n$-module) 
is isomorphic to $K \otimes_A \VC_\LG$ (resp. $\QM \otimes_\ZM \widehat{\VC}_\LG$) 
for some left cell $\LG$.

\end{itemize}
\end{prop}

\begin{cor}\label{constructible}
An $\HC_n^K$-module is constructible if and only if it is isomorphic 
to $K \otimes_A \VC_\LG$ for some left cell $\LG$ of $W_n$.
\end{cor}

\begin{proof}
This follows from \cite[Prop.~5.2]{Ge1} (for the definition of a constructible 
module for $\HC_n^K$, the reader may refer to \cite[\S 22]{Lusztig02}).
\end{proof}

We conclude by giving an explicit description of the irreducible 
character of $\HC_n^K$ (or $W_n$) defined by a left cell. We need some notation. 

We denote by $\e_n$ the sign character of $W_n$ (it is defined 
by $\e_n(t)=\e_n(s_1)=\dots=\e_n(s_{n-1})=-1$). If $\l$ is a partition 
of $n$, we denote by $1_\l$ (resp. $\e_\l$) the restriction to $\SG_\l$ 
of the trivial character (resp. of $\e_n$) to $\SG_\l$. We then 
define $\ch_\l \in \Irr \SG_n$ to be the unique common irreducible component 
of $\Ind_{\SG_\l}^{\SG_n} 1_\l$ and $\Ind_{\SG_{\l^*}}^{\SG_n} \e_{\l^*}$. 
We denote by $\ch_\l^+$ the irreducible character of $W_n$ obtained 
by composition of $\ch_\l$ with the surjective morphism $W_n \to \SG_n$. 

If $0 \leq l \leq n$, we denote by $\th_{l,n-l}$ the linear character 
of the normal subgroup $N_n=<t_1,t_2,\dots,t_n>$ of $W_n$ such that 
$$\th_{l,n-l}(t_i)=\left\{\begin{array}{rl}
                          1 & \text{if}~i \leq l \\
			  -1 & \text{if}~i \geq l+1.
			  \end{array}\right.$$
The stabilizer of $\th_{l,n-l}$ in $\SG_n$ is $\SG_{l,n-l}$. 
We denote by $\thet_{l,n-l}$ the linear character of $\SG_{l,n-l}.N_n$ 
whose restriction to $N_n$ is $\th_{l,n-l}$ and which is trivial 
on $\SG_{l,n-l}$. Now, if $\l$ is a partition of $n$, we denote 
by $\ch_\l^-$ the irreducible character $\ch_\l^+.\thet_{0,n}$ of $W_n$.

We say that an irreducible character $\ch$ of $W_n$ {\it has weight $l$} 
if $\th_{l,n-l}$ occurs in $\Res_{N_n}^{W_n} \ch$. 
By Clifford theory, the weight of an irreducible character of $W_n$ 
is uniquely determined.

If $(\l,\m)$ is a bipartition of $n$, we denote by $\ch_\l^+ \boxtimes \ch_\m^-$ 
the irreducible character of $\SG_{|\l|,|\m|}.N_n \simeq W_{|\l|} \times W_{|\m|}$ 
obtained by the (external) tensor product of $\ch_\l^+ \in \Irr W_{|\l|}$ and 
$\ch_\m^- \in \Irr W_{|\m|}$. We then set
$$\ch_{\l,\m}^W=\Ind_{\SG_{|\l|,|\m|}.N_n}^{W_n} 
\Bigl(\ch_\l^+ \boxtimes \ch_\m^-\Bigr).$$
By Clifford theory, $\ch_{\l,\m}^W$ is an irreducible character of $W_n$ 
and the map $(\l,\m) \mapsto \ch_{\l,\m}^W$ is a bijection between the 
set of bipartitions of $n$ and $\Irr W_n$. Note that $\ch_{\l,\m}^W$ 
has weight $|\l|$. If $\l$ is a partition of $n$, then $\ch_{\l,0}^W=\ch_\l^+$ 
and $\ch_{0,\l}^W=\ch_\l^-$. 

Let us now talk about irreducible characters of Hecke algebras. 
We denote by $\HC(\SG_n)$ the sub-$A$-algebra of $\HC_n$ 
generated by $T_{s_1}$, $T_{s_2}$, \dots, $T_{s_{n-1}}$. We set 
$\HC^K(\SG_n) = K \otimes_A \HC(\SG_n)$.  
If $\l$ is a partition of $n$, we denote by $\ch_\l^\HC$ the unique irreducible 
character of $\HC^K(\SG_n)$ such that $\widehat{\ch_\l^\HC(T_w)}=\ch_\l(w)$ 
for every $w \in \SG_n$. If $(\l,\m)$ is a bipartition of $n$, we denote 
by $\ch_{\l,\m}^\HC$ the unique irreducible character of $\HC_n$ such that 
$\widehat{\ch_{\l,\m}^\HC(T_w)}=\ch_{\l,\m}^W(w)$ for every $w \in W_n$. 

Finally, if $\LG$ is a left cell, we define the {\it shape} of $\LG$ to be 
the bipartition $(\l,\m)$ such that $B_n^+(w)$ has shape $\l$ and 
$B^-_n(w)$ has shape $\m$ for every $w \in \LG$. 

\begin{prop}\label{irr cell}
Let $\LG$ be a left cell of shape $(\l,\m)$. Then the irreducible 
character of $\HC_n^K$ afforded by $K \otimes_A \VC_\LG$ is 
$\ch_{\m,\l^*}^\HC$.
\end{prop}

\begin{proof}
By Proposition \ref{left droite}, we may (and we will) assume that 
$b_w=1$ for every $ w\in \LG$. We set $l=|\m|$, so that $|\l|=n-l$. 
We denote by $\ch$ the irreducible character of $W_n$ afforded 
by $\QM \otimes_\ZM \widehat{\VC}_\LG$. 
We denote by $\HC_{l,n-l}$ the sub-$A$-algebra of $\HC_n$ generated 
by $(T_s)_{s \in S_{l,n-l}}$. It is a Hecke algebra for the Weyl 
group $W_{l,n-l} \simeq W_l \times \SG_{n-l}$. We also set 
$\HC_{l,n-l}^K=K \otimes_A \HC_{l,n-l}$. 

Let $\LG'$ (resp. $\LG''$) be the left cell of $\SG_l$ (resp. $\SG_{[l+1,n]}$) 
such that $\s_w^\prime \in \LG'$ (resp. $\s_w^{\prime\prime} \in \LG''$) 
for every $w \in \LG$ (see Theorem \ref{main theo}). If 
$w \in \LG$, we denote by $c_w$ the image of $C_w \in \IC_{\leq \LG}$ 
in $\VC_\LG=\IC_{\leq \LG}/\IC_{< \LG}$. 

Then $\LG'$ has shape $\m$ and $\LG''$ has shape $\l$. We set 
$$\MG=\{w \in \LG~|~a_w=1\}.$$
We denote by $M_\LG$ the $A$-submodule of $\VC_\LG$ generated 
by $(c_w)_{w \in \MG}$. If $w \in \MG$, then $C_w = C_{w_l \s_w^\prime} 
C_{\s_w^{\prime\prime}}$. 
This shows that $M_\LG$ is a sub-$\HC_{l,n-l}$-module isomorphic to 
$\VC_{w_l \LG'} \boxtimes_A \VC_{\LG''}$ through the canonical isomorphism 
$\HC_{l,n-l} \simeq \HC_l \otimes_A \HC(\SG_{[l+1,n]})$ 
(note that $\LG'$ is also a left cell of $W_l$). 

But, the character of $W_l$ afforded by $\QM \otimes_\ZM 
\widehat{\VC}_{w_l \LG'}$ 
is equal to the product of $\e_l$ with the character afforded by 
$\QM \otimes_\ZM \widehat{\VC}_{\LG'}$. 
Since $C_t C_w = C_{tw}$ for every $w \in \SG_n$, 
we get that $C_t=T_t+Q^{-1}$ acts as $0$ on $\VC_{\LG'}$. So the character of $W_n$ 
afforded by $\QM \otimes_\ZM \widehat{\VC}_{\LG'}$ is $\ch_{\m^*}^-$. In other words, 
the character of $W_l$ afforded by $\QM \otimes_\ZM \widehat{\VC}_{w_l \LG'}$ 
is $\e_l.\ch_{\m^{*}}^{-}=\ch_\m^+$. Therefore,
$$\langle \Res_{W_{l,n-l}}^{W_n} \ch , 
\ch_\m^+ \boxtimes \ch_{\l^*} \rangle_{W_{l,n-l}}\not= 0.\leqno{(*)}$$
This shows that the weight of $\ch$ is greater than or equal to $l$. 

But, if we use the same argument for the left cell $w_n \LG$ of $W_n$, we get 
that the weight of $\e_n.\ch$ is greater than or equal to $n-l$. Therefore, 
the weight of $\ch$ is equal to $l$. By $(*)$, $\ch$ is an irreducible constituent 
of weight $l$ of 
$$\Ind_{\SG_{l,n-l}.N_n}^{W_n} 
\Bigl(\ch_\m^+ \boxtimes (\Ind_{\SG_{n-l}}^{W_{n-l}} \ch_{\l^*})\Bigr).$$
But, $\ch_{\l^*}^-$ is the unique irreducible constituent of 
$\Ind_{\SG_{n-l}}^{W_{n-l}} \ch_{\l^*}$ 
of weight $0$. Therefore, $\ch$ is an irreducible constituent of 
$$\Ind_{\SG_{l,n-l}.N_n}^{W_n} \Bigl(\ch_\m^+ \boxtimes \ch_{\l^*}^-\Bigr).$$
Consequently, $\ch=\ch_{\m,\l^*}^W$.
\end{proof}
\end{what}

\begin{what}\label{final !}{\bf Final remarks.} 
We conclude by mentioning three further developments, that will make the 
object of a forthcoming paper.

\medskip

{\bf 1.} The generalized Robinson-Schensted correspondence described in 
\S \ref{sec2} can be extended for the complex reflection groups 
$G(e,1,n) = (\ZM/e\ZM) \wr \SG_n$ (see \cite{iancu}). RS-cell 
representations of the corresponding Ariki-Koike algebra 
can be constructed as in Proposition \ref{irr cell}. 

\medskip

{\bf 2.} In the asymptotic case that we have studied here, the Kazhdan-Lusztig 
basis proves to be a {\it cellular basis} for $\HC_n$ in the sense 
of Graham-Lehrer \cite[Definition 1.1]{GraLe}. We intend to investigate 
the links between the Kazhdan-Lusztig basis and the cellular basis 
constructed by Graham and Lehrer \cite[\S 5]{GraLe}. 

\medskip

{\bf 3.} We expect that our explicit results will help to prove 
Lusztig's conjectures (P4), (P9), (P10) and (P11) \cite[\S 14]{Lusztig02} 
in this asymptotic context.
\end{what}


\end{document}